\documentclass[11pt, twoside]{article}
\usepackage{latexsym}
\usepackage{amsmath}
\usepackage{amssymb}
\usepackage[all]{xy}
\usepackage{amsfonts}
\usepackage{verbatim}
\usepackage{amsthm}
\usepackage{mathrsfs}
\usepackage{epsfig}
\usepackage{xy}
\usepackage{tensor}
\usepackage{array}
\usepackage{stmaryrd}
\usepackage{graphicx,color}
\usepackage{xcolor}
\usepackage[colorlinks=true,linkcolor=blue,citecolor=blue]{hyperref}
\usepackage{tikz}
\usetikzlibrary{arrows,calc}
\usepackage{etex}
\usepackage{mathdots}
\usepackage{float}
\usepackage{graphics}
\usepackage{pdflscape}
\usepackage{extarrows}
\usepackage{relsize}
\usepackage{anysize,hyperref}
\input xypic
\xyoption{all}
\usepackage{bm}
\usepackage[perpage,symbol]{footmisc}
\usepackage{setspace}
\SelectTips{cm}{10}

\usepackage{enumitem}
\setenumerate[1]{itemsep=0pt,partopsep=0pt,parsep=\parskip,topsep=5pt}
\setitemize[1]{itemsep=0pt,partopsep=0pt,parsep=\parskip,topsep=5pt}
\setdescription{itemsep=0pt,partopsep=0pt,parsep=\parskip,topsep=5pt}
\renewcommand{\cal}{\mathcal}
\def\A{\mathscr{A}}

\def\M{\mathscr{M}}

\def\Id{\mathrm{Id}}

\def\dr{\ar@{->}[r]}
\def\Im{\mbox{\rm Im}\,}\def\Ker{\mbox{\rm Ker}\,}

\def\X{\mathscr{X}}

\def\add{\mbox{add}}
\def\Ext{\mbox{Ext}}
\def\Hom{\mbox{Hom}}

\def\modcat{\mbox{mod\text{-}}}

\begin{document}
\baselineskip=15pt
\title{\Large{\bf  Wakamatsu tilting subcategories, weak support $\tau$-tilting subcategories and recollements  \footnotetext{$^\ast$Corresponding author. }}}
\medskip
\author{ Yongduo Wang, Hongyang Luo, Yu-Zhe Liu, Jian He$^\ast$ and Dejun Wu }

\date{}

\maketitle
\def\blue{\color{blue}}
\def\red{\color{red}}

\newtheorem{theorem}{Theorem}[section]
\newtheorem{lemma}[theorem]{Lemma}
\newtheorem{corollary}[theorem]{Corollary}
\newtheorem{proposition}[theorem]{Proposition}
\newtheorem{conjecture}{Conjecture}
\theoremstyle{definition}
\newtheorem{definition}[theorem]{Definition}
\newtheorem{question}[theorem]{Question}
\newtheorem{notation}[theorem]{Notation}
\newtheorem{remark}[theorem]{Remark}
\newtheorem{remark*}[]{Remark}
\newtheorem{example}[theorem]{Example}
\newtheorem{example*}[]{Example}
\newtheorem{condition}[theorem]{\bf Condition}
\newtheorem{construction}[theorem]{Construction}
\newtheorem{construction*}[]{Construction}

\newtheorem{assumption}[theorem]{Assumption}
\newtheorem{assumption*}[]{Assumption}

\baselineskip=17pt
\parindent=0.5cm

\begin{abstract}
\baselineskip=16pt
In this article, we prove that if ($\mathcal{A}$, $\mathcal{B}$, $\mathcal{C}$) is a recollement of abelian categories, then Wakamatsu tilting (resp. weak support $\tau$-tilting) subcategories in $\mathcal{A}$ and $\mathcal{C}$ can induce Wakamatsu tilting (resp. weak support $\tau$-tilting) subcategories in  $\mathcal{B}$, and the converses hold under natural assumptions. As an application, we mainly consider the relationship of $\tau$-cotorsion torsion triples in  ($\mathcal{A}$, $\mathcal{B}$, $\mathcal{C}$).        \\[0.5cm]
\medskip
\end{abstract}
\pagestyle{myheadings}
\markboth{\rightline {\scriptsize Y. Wang, H. Luo, Y. Liu, J. He and D. Wu  }}
         {\leftline{\scriptsize   Wakamatsu tilting subcategories, weak support $\tau$-tilting subcategories and recollements}}

\section{Introduction}
The recollements of abelian categories first appeared in the construction of the
category of perverse sheaves on a singular space in \cite{BBD}.  It has been applied to many aspects of algebra, for example, representation theory, ring theory, geometry etc.

Tilting theory plays an important role in the representation theory of Artin algebras, see \cite{AR, H}. The notion
of the classical tilting modules over an Artin algebra was introduced by Brenner and Butler in \cite{BB}. Since then, tilting
modules have been investigated by many authors. In some sense, cotilting theory is the dual of tilting theory. A further generalization of tilting modules to modules of possibly infinite projective dimension was made by
Wakamatsu \cite{W}, which is now called Wakamatsu tilting module. Note that Wakamatsu tilting modules are
common generalizations of tilting modules and cotilting modules. Wakamatsu tilting subcategory is a certain
categorical analogue of Wakamatsu tilting module defined in \cite{E}. Suppose that $\mathcal{B}$ admits a recollement relative to abelian categories $\mathcal{A}$ and $\mathcal{C}$. Our first
main result describes how to glue together Wakamatsu tilting subcategory in $\mathcal{A}$ and $\mathcal{C}$, to
obtain a Wakamatsu tilting subcategory of $\mathcal{B}$, see Theorem \ref{main1}. In the reverse direction, we
give sufficient conditions on a Wakamatsu tilting subcategory of $\mathcal{B}$, relative to the functors
involved in the recollement, to Wakamatsu tilting subcategories in $\mathcal{A}$ and $\mathcal{C}$, see Theorem \ref{main2}.

Iyama and Reiten introduced $\tau$-tilting theory \cite{AIR}, which is a generalization
of classical tilting theory. Due to the effectiveness of $\tau$-tilting theory for the study of the categories of finitely presented
modules, many authors have introduced theories generalizing $\tau$-tilting theory, and its
dual, to other contexts, see \cite{AMV, AST, BP, IJY, LZ}. Particularly, Asadollahi, Sadeghi and Treffinger \cite{AST} studied $\tau$-tilting theory in arbitrary abelian categories with
enough projective objects. Meanwhile, they introduced the concept of weak support $\tau$-tilting (resp. support $\tau$-tilting) subcategories in an abelian category. Our second
main result describes how to glue together weak support $\tau$-tilting subcategories in $\mathcal{A}$ and $\mathcal{C}$, to
obtain a weak support $\tau$-tilting subcategory of $\mathcal{B}$, see Theorem \ref{main3}. In the reverse direction, we
give sufficient conditions on a weak support $\tau$-tilting subcategory of $\mathcal{B}$, relative to the functors
involved in the recollement, to weak support $\tau$-tilting subcategories in $\mathcal{A}$ and $\mathcal{C}$, see Theorem \ref{main4}. As Corollaries of Theorem \ref{main3} and \ref{main4}, we also glue support $\tau$-tilting subcategory in a recollement ($\mathcal{A}$, $\mathcal{B},\mathcal{C}$) of abelian categories, see Corollary  \ref{xx} and \ref{bb}. Note that Asadollahi, Sadeghi and Treffinger \cite{AST} established a bijection between $\tau\textrm{-cotorsion torsion triples}$ and $\textrm{support}~\tau\textrm{-tilting subcategories}$. As an application, we finally consider the relationship of $\tau$-cotorsion torsion triples in  ($\mathcal{A}$, $\mathcal{B}$, $\mathcal{C}$), see Proposition \ref{p1} and \ref{p2}.

This article is organized as follows. In Section 2, we give some terminologies
and some preliminary results. In Section 3, we prove our first main
result. In Section 4, we prove our second main
result. In Section 5, we give some applications of our main results. In Section 6, we give some examples to explain our main results.

\section{Preliminaries}
In this paper, all subcategories are  full subcategories, closed under isomorphisms and direct sums. First, let us recall the concept of a recollement of abelian categories from \cite{FP}.

\begin{definition}\label{recollement}{\rm \cite{FP}}
Let $\mathcal{A}$, $\mathcal{B}$ and $\mathcal{C}$ be three abelian categories. A \emph{recollement} of $\mathcal{B}$ relative to
$\mathcal{A}$ and $\mathcal{C}$, denoted by ($\mathcal{A}$, $\mathcal{B}$, $\mathcal{C}$), is a diagram
\begin{equation}\label{recolle}
  \xymatrix{\mathcal{A}\ar[rr]|{i_{*}}&&\ar@/_1pc/[ll]|{i^{*}}\ar@/^1pc/[ll]|{i^{!}}\mathcal{B}
\ar[rr]|{j^{\ast}}&&\ar@/_1pc/[ll]|{j_{!}}\ar@/^1pc/[ll]|{j_{\ast}}\mathcal{C}}
\end{equation}
given by two exact functors $i_{*},j^{\ast}$, two right exact functors $i^{\ast}$, $j_!$ and two left exact functors $i^{!}$, $j_\ast$, which satisfies the following conditions:
\begin{itemize}
  \item [(R1)] $(i^{*}, i_{\ast}, i^{!})$ and $(j_!, j^\ast, j_\ast)$ are adjoint triples.
  \item [(R2)] $\Im i_{\ast}=\Ker j^{\ast}$.
  \item [(R3)] $i_\ast$, $j_!$ and $j_\ast$ are fully faithful.
\end{itemize}
\end{definition}

Next, we collect some properties of recollements (see \cite{FP, P, PSS, PV}), which are very useful in the sequel.

\begin{lemma}\label{7} \rm Let ($\mathcal{A}$, $\mathcal{B}$, $\mathcal{C}$) be a recollement of abelian categories as \rm{(\ref{recolle})}.

$(1)$ All the natural transformations
$$i^{\ast}i_{\ast}\Rightarrow\Id_\mathcal{A},~\Id_\mathcal{A}\Rightarrow i^{!}i_{\ast},~\Id_\mathcal{C}\Rightarrow j^{\ast}j_{!},~j^{\ast}j_{\ast}\Rightarrow\Id_\mathcal{C}$$
are natural isomorphisms.

$(2)$ $i^{\ast}j_!=0$ and $i^{!}j_\ast=0$.

$(3)$ If $i^{\ast}$ (resp. $i^{!}$) is  exact, then ${i^!}{j_!} = 0$ (resp. ${i^*}{j_*} = 0$).

$(4)$ If $i^{\ast}$ (resp. $i^{!}$) is exact, then $j_{!}$ (resp. $j_{\ast}$) is  exact.

$(5)$ For any $B \in {\rm{{\cal B}}}$, if ${i^*}$ is exact, then there is an exact sequence \begin{equation}\xymatrix{0\ar[r]&{j_!}{j^*}(B)\ar[r]&B\ar[r]&{i_*}{i^*}(B)\ar[r] &0}.\label{exa1}\end{equation}

$(5')$ For any $B \in {\rm{{\cal B}}}$, if ${i^!}$ is exact, then there is an exact sequence \begin{equation}\xymatrix{0\ar[r]&{i_*}{i^!}(B)\ar[r]&B\ar[r]&{j_*}{j^*}(B)\ar[r] &0}.\label{exa2}\end{equation}

\end{lemma}

\begin{remark}For a recollement of abelian categories as \rm{(\ref{recolle})}, if all the functors are exact, then we call this situation a \emph{good recollement}. In order to achieve the situation, we actually need to verify the exactness only for two functors $i^{\ast}$ and $i^{!}$ by Lemma \ref{7}.
\end{remark}

It is well known that let $f:\mathcal{A}\rightarrow\mathcal{B}$ and $g: \mathcal{B}\rightarrow\mathcal{A}$ be two functors between two abelian categories $\mathcal{A}$ and $\mathcal{B}$, if $(g, f)$ is a adjoint pair and $f$ is exact, then $g$ preserves projective objects. In fact, let $P$ be a projective object in $\mathcal{B}$ and $0\longrightarrow X\stackrel{}\longrightarrow Y\stackrel{}\longrightarrow Z\longrightarrow0$ an exact sequence in $\mathcal{A}$. And there is a commutative diagram
 $$\xymatrix@C=1.2cm{
0\ar[r] &\mathcal{A}(gP,X)\ar[r]\ar[d]^{\cong} & \mathcal{A}(gP,Y) \ar[r]\ar[d]^{\cong} & \mathcal{A}(gP,Z) \ar[r]^{}\ar[d]^{\cong}& 0\\
0\ar[r] &\mathcal{B}(P,fX)\ar[r] & \mathcal{B}(P,fY) \ar[r]& \mathcal{B}(P,fZ)\ar[r]^{}& 0,
}$$
where the second row is an exact sequence. So the first row is also an exact sequence, i.e., $gP$ is a projective object in $\mathcal{A}$. Moreover, we have the following Lemma \ref{4}.
\begin{lemma}\label{4}\rm{\cite[Lemma 3.10]{L}} If $(g, f)$ is a adjoint pair as above, $\mathcal{B}$ has enough projectives,
both $f$ and $g$ are exact, then $\mathrm{Ext}_\mathcal{A}^n(g(B),A)\cong \mathrm{Ext}_\mathcal{B}^n(B, f(A))$, for all $A\in \mathcal{A}$, $B\in \mathcal{B}$ and any positive integer $n$.
\end{lemma}

\begin{remark} \label{Yu}For a good recollement of abelian categories as \rm{(\ref{recolle})}, if $\mathcal{A}$ and $\mathcal{B}$ have enough projectives, then $$\mathrm{Ext}_\mathcal{B}^1(i_*i^*(B),j_!j^*(B))\cong \mathrm{Ext}_\mathcal{A}^1(i^*(B),i^!j_!j^*(B))=0$$ and $$\mathrm{Ext}_\mathcal{B}^1(j_*j^*(B),i_*i^!(B))\cong \mathrm{Ext}_\mathcal{A}^1(i^*j_*j^*(B),i^!(B))=0$$ by Lemmas \ref{7} and \ref{4}. That is to say, the exact sequences (\ref{exa1}) and (\ref{exa2}) are split in Lemma \ref{7}.
\end{remark}

\begin{lemma}\label{horse}\rm{\cite[Lemma 2.1]{PHZ}}
\rm Let $\mathcal{A}$ be an abelian category and $0\longrightarrow X\stackrel{f}\longrightarrow Y\stackrel{g}\longrightarrow Z\longrightarrow0$ an exact sequence in $\mathcal{A}$.

\rm(1) Assume that $$\xymatrix{0\ar[r]&X\ar[r]^{a}&X_{0}\ar[r]^{a_{0}}&X_{-1}\ar[r]^{a_{-1}}&\cdots}$$
is a complex and $$\xymatrix{0\ar[r]&Z\ar[r]^{b}&Z_{0}\ar[r]^{b_{0}}&Z_{-1}\ar[r]^{b_{-1}}&\cdots}$$ an exact sequence in $\A$.
If $\Ext_{\mathcal{A}}^{1}(\Ker\,b_{i},~X_{i})=0$ for any $i\leq 0$, then there exist diagram with exact rows
\begin{center}$\xymatrix{&0\ar[d]&0\ar[d]&0\ar[d]&\\
                         0\ar[r]&X\ar[r]^{f}\ar[d]^{a}&Y\ar[r]^{g}\ar[d]&Z\ar[r]\ar[d]^{b}&0\\
                         0\ar[r]&X_{0}\ar[r]\ar[d]^{a_{0}}&Z_{0}\oplus X_{0}\ar[r]\ar[d]&Z_{0}\ar[r]\ar[d]^{b_{0}}&0\\
                         0\ar[r]&X_{-1}\ar[r]\ar[d]^{a_{-1}}&Z_{-1}\oplus X_{-1}\ar[r]\ar[d]&Z_{-1}\ar[r]\ar[d]^{b_{-1}}&0\\
                         &\vdots & \vdots & \vdots &}$\end{center}
commutes. Moreover, the middle column is exact if and only if the left column is exact.

\rm(2) Assume that $$\xymatrix{\cdots \ar[r]& X'_{1}\ar[r]^{a'_{1}}& X'_{0}\ar[r]^{a'_{0}}&X\ar[r]&0}$$ be an exact sequence
and $$\xymatrix{\cdots \ar[r]& Z'_{1}\ar[r]^{b'_{1}}&}\xymatrix{Z'_{0}\ar[r]^{b'_{0}}&Z\ar[r]&0}$$ a complex in $\mathcal{A}$.
If $\Ext_{\mathcal{A}}^{1}(Z'_{i},\Im\,a'_{i})=0$ for any $i\geq 0$, then there exists diagram with exact rows
 \begin{center}$\xymatrix{&\vdots\ar[d] & \vdots\ar[d] & \vdots\ar[d] &\\
 0\ar[r]&X'_{1}\ar[r]\ar[d]^{a'_{1}}&Z'_{1}\oplus X'_{1}\ar[r]\ar[d]&Z'_{1}\ar[r]\ar[d]^{b'_{1}}&0\\
 0\ar[r]&X'_{0}\ar[r]\ar[d]^{a'_{0}}&Z'_{0}\oplus X'_{0}\ar[r]\ar[d]&Z'_{0}\ar[r]\ar[d]^{b'_{0}}&0\\
 0\ar[r]&{X}\ar[r]^{f}\ar[d]&{Y}\ar[d]\ar[r]^{g}&{Z}\ar[d]\ar[r]&0\\
 &0 & 0 & 0 &}$\end{center}
 commutes. Moreover, the middle column is exact if and only if the right column is exact.
\end{lemma}

\section{Gluing Wakamatsu tilting subcategories }
Suppose that $\mathcal{C}$ is a category, for any $A,B\in\mathcal{C}$, we simply denote ${{\rm Hom}_{\mathcal{C}}(A,B)}$ as ${\mathcal{C}(A,B)}$ in this paper.
From now until the end of the article, all of the abelian categories we will consider have enough projective objects and injective objects.

Suppose that $\mathcal{A}$ is an abelian category and $\mathcal{D}$ is a subcategory of $\mathcal{A}$.

Set
$${}^ \bot \mathcal{D}: = \{M \in \mathcal{A}~|~ \mathrm{Ext}_\mathcal{A}^i( {M,N} ) = 0,\forall i \ge 1,\forall N \in \mathcal{D} \}.$$
The subcategory $\mathcal{D}$ is said to be self-orthogonal if $\mathcal{D} \subseteq {}^ \bot \mathcal{D}$.

Suppose $\mathcal{W}$ is a subcategory of $\mathcal{A}$, the full subcategory $X_{\mathcal{W}}$ of $\mathcal{A}$ is defined as follow:
{\small $$
X_{\mathcal{W}}:=\biggl\{~M\in{}^ \bot \mathcal{W}\bigg|\begin{array}{ll}
 \mbox{there is an exact sequence }0\xrightarrow{}M\xrightarrow{d_0} W_0\xrightarrow{d_1} W_1\xrightarrow{} \cdots\\
\textrm{for}~i \geq 0~\textrm{with}~{W_i} \in \mathcal{W},\mathrm{Im}{\kern 1pt} {d_i} \in {}^ \bot {\rm{{\cal W}}}
 \end{array}\biggl\}.$$}
We first recall the concept of Wakamatsu tilting subcategory in an abelian category.
\begin{definition}\label{Wakamatsu}{\rm \cite[Definition 3.1]{E}}
We say that $\mathcal{W}$ is a Wakamatsu tilting subcategory of $\mathcal{A}$ if it satisfies the following conditions:
\begin{itemize}
  \item [(1)]  $\mathcal{W}$ is self-orthogonal.
  \item [(2)]  $X_\mathcal{W}$ contains all projective objects in $\mathcal{A}$.
\end{itemize}
\end{definition}


Our first main result is the following.
\begin{theorem}\label{main1}\rm
Let ($\mathcal{A}$, $\mathcal{B}$, $\mathcal{C}$) be a good recollement of abelian categories as \rm{(\ref{recolle})}, ${\rm{{\cal X}}}'$ and ${\rm{{\cal X}}}''$ are Wakamatsu tilting subcategories of $\mathcal{A}$ and $\mathcal{C}$. Define $$\mathcal{X}= \{ {X \in \mathcal{B}~|~{i^!}( X ) \in \mathcal{X}',{j^*}( X ) \in \mathcal{X}''} \}.$$ Then $\mathcal{X}$ is a Wakamatsu tilting subcategory of $\mathcal{B}$.
\end{theorem}

\begin{proof}
We first claim that $\mathcal{X}$ is a self-orthogonal subcategory of ${\rm{{\cal B}}}$. Indeed, let $M$, $N \in {\rm{{\cal X}}}$, there are two exact sequences
$$\xymatrix{0\ar[r]&{i_*}{i^!}(M)\ar[r]&M\ar[r]&{j_*}{j^*}(M)\ar[r] &0}$$
and
$$\xymatrix{0\ar[r]&{i_*}{i^!}(N)\ar[r]&N\ar[r]&{j_*}{j^*}(N)\ar[r] &0}$$
by Lemma \ref{7}. Hence we have the following two exact sequences
$$\xymatrix{{\rm Ext}_\mathcal{B}^n(j_*j^*(M),N)\ar[r]&{\rm Ext}_\mathcal{B}^n(M,N)\ar[r]&{\rm Ext}_\mathcal{B}^n(i_*i^!(M),N)}$$
and$$\xymatrix{{\rm Ext}_\mathcal{B}^n(j_*j^*(M),i_*i^!(N))\ar[r]&{\rm Ext}_\mathcal{B}^n(j_*j^*(M),N)\ar[r]&{\rm Ext}_\mathcal{B}^n(j_*j^*(M),j_*j^*(N))}.$$
Note that ${i^!}( M ) \in \mathcal{X}',{j^*}( M) \in \mathcal{X}'',{i^!}( N ) \in \mathcal{X}',{j^*}( N) \in \mathcal{X}''$.
From Lemma \ref{4}, we have $$0 = \mathrm{Ext}_\mathcal{A}^n( {i^!}( M ),{i^!}( N ) ) \cong \mathrm{Ext}_\mathcal{B}^n( {i_*}{i^!}( M ),N ),$$
$$0 = \mathrm{Ext}_\mathcal{A}^n( 0,{i^!}( N ) ) \cong \mathrm{Ext}_\mathcal{A}^n( {i^*}{j_*}{j^*}( M ),{i^!}( N ) ) \cong \mathrm{Ext}_\mathcal{B}^n( {j_*}{j^*}( M ),{i_*}{i^!}( N ) )$$ and
$$0 = \mathrm{Ext}_\mathcal{C}^n( {j^*}( M ),{j^*}( N ) ) \cong \mathrm{Ext}_\mathcal{C}^n( {j^*}{j_*}{j^*}( M ),{j^*}( N ) ) \cong \mathrm{Ext}_\mathcal{B}^n( {j_*}{j^*}( M ),{j_*}{j^*}( N ))$$ hold.
So $\mathrm{Ext}_\mathcal{B}^n( {j_*}{j^*}( M ),N ) = 0$, then $\mathrm{Ext}_\mathcal{B}^n( M,N ) = 0$. This shows that $\mathcal{X}$ is self-orthogonal subcategory of $\mathcal{B}$. Next, we prove that $X_\mathcal{X}$ contains all projective objects in $\mathcal{B}$.

Let $P$ be a projective object of ${\rm{{\cal B}}}$, we know that ${i^*}\left( P \right)$ and ${j^*}\left( P \right)$ are projective objects of $\mathcal{A}$ and $\mathcal{C}$, respectively. So
${i^*}\left( P \right) \in {X_{{\rm{{\cal X}}}'}}$, ${j^*}\left( P \right) \in {X_{{\rm{{\cal X}}}''}}$. Then there are two exact sequences
{$$\xymatrix{0\ar[r]&{i^*}(P)\ar[r]^{d_0}&{T_0}\ar[r]^{d_1}&{T_1}\ar[r] &\cdots}$$ }and
$$\xymatrix{0\ar[r]&{j^*}(P)\ar[r]^{q_0}&{Q_0}\ar[r]^{q_1}&{Q_1}\ar[r] &\cdots}$$
with ${T_i} \in {\rm{{\cal X}}}'$, ${Q_i} \in {\rm{{\cal X}}}''$, ${\mathop{\rm Im}\nolimits} {\kern 1pt} {d_i} \in {}^ \bot {\rm{{\cal X}}}'$ and
${\mathop{\rm Im}\nolimits} {\kern 1pt} {q_i} \in {}^ \bot {\rm{{\cal X}}}''$ for any $i \ge 0$. Therefore, there are two exact sequences
$$\xymatrix{0\ar[r]&{i_*}{i^*}(P)\ar[r]^{{i_*}({d_0})}&{i_*}({T_0})\ar[r]^{{i_*}({d_1})}&{i_*}({T_1})\ar[r] &\cdots}$$ and
$$\xymatrix{0\ar[r]&{j_!}{j^*}(P)\ar[r]^{{j_!}({q_0})}&{j_!}({Q_0})\ar[r]^{{j_!}({q_1})}&{j_!}({Q_1})\ar[r] &\cdots}.$$
By Lemma \ref{7}, there exists an exact sequence
$$\xymatrix{0\ar[r]&{j_!}{j^*}(P)\ar[r]&P\ar[r]&{i_*}{i^*}(P)\ar[r] &0}.$$
Note that
$${\mathop{\rm Ext}\nolimits} _{\rm{{\cal B}}}^1(\Ker{i_*}\left( {{d_1}} \right),\;{j_!}({Q_0})) \cong {\mathop{\rm Ext}\nolimits} _{\rm{{\cal B}}}^1\left( {{i_*}{i^*}\left( P \right),{j_!}({Q_0})} \right) \cong {\mathop{\rm Ext}\nolimits} _{\rm{{\cal A}}}^1( {{i^*}\left( P \right),{i^!}{j_!}({Q_0})} ) = 0.$$
Then there is a commutative diagram with exact rows and columns by Lemma \ref{horse}
\begin{center}$\xymatrix{&0\ar[d]&0\ar[d]&0\ar[d]&\\
                         0\ar[r]&j_!j^*(P)\ar[r]\ar[d]&P\ar[r]\ar[d]^{\varepsilon_0}&i_*i^*(P)\ar[r]\ar[d]&0\\
                         0\ar[r]&j_!(Q_0)\ar[r]\ar[d]&i_*(T_0)\oplus j_!(Q_0)\ar[r]\ar[d]&i_*(T_0)\ar[r]\ar[d]&0\\
                         0\ar[r]&j_!(\mathrm{Im}q_1)\ar[r]\ar[d]&C_0\ar[r]\ar[d]&i_*(\mathrm{Im}d_1)
                         \ar[r]\ar[d]&0\\
                         & 0 & 0 & 0 }$\end{center}
where ${C_0} = {\mathop{\rm Coker}\nolimits} \varepsilon_0 $. Hence, there is an exact sequence
$$\xymatrix{0\ar[r]&P\ar[r]&{i_*}({T_0})\oplus {j_!}({Q_0})\ar[r]&C_0\ar[r] &0}.$$
Since ${\rm{{\cal X}}}'$ and ${\rm{{\cal X}}}''$ are closed to direct sums and isomorphisms, and
$${i^!}\left( {{i_*}({T_0}) \oplus {j_!}({Q_0})} \right) \cong {i^!}{i_*}\left( {{T_0}} \right) \cong {T_0},$$
$${j^*}\left( {{i_*}({T_0}) \oplus {j_!}({Q_0})} \right) \cong {j^*}{j_!}\left( {{Q_0}} \right) \cong {Q_0},$$
then $\left( {{i_*}({T_0}) \oplus {j_!}({Q_0})} \right) \in {\rm{{\cal X}}}$.
Finally, we show that ${C_0} \in {}^ \bot {\rm{{\cal X}}}$. Let $N \in {\rm{{\cal X}}}$, note that the exact sequence
$$\xymatrix{0\ar[r]&{j_!}({\mathop{\rm Im}\nolimits} {q_1})\ar[r]&{C_0}\ar[r]&{i_*}({\mathop{\rm Im}\nolimits} {d_1})\ar[r] &0},$$
there is an exact sequence
$$\xymatrix{\mathrm{Ext}_\mathcal{B}^n(i_*(\mathrm{Im}d_1),N)\ar[r]&\mathrm{Ext}_\mathcal{B}^n(C_0,N)\ar[r]&\mathrm{Ext}_\mathcal{B}^n(j_!(\mathrm{Im}q_1),N)}.$$
While $${\mathop{\rm Ext}\nolimits} _\mathcal{B}^n\left( {{j_!}({\mathop{\rm Im}\nolimits} {q_1}),N} \right) \cong {\mathop{\rm Ext}\nolimits} _{\rm{{\cal C}}}^n\left( {{\mathop{\rm Im}\nolimits} {q_1},{j^*}\left( N \right)} \right) = 0$$ and
$${\mathop{\rm Ext}\nolimits} _\mathcal{B}^n\left( {{i_*}({\mathop{\rm Im}\nolimits} {d_1}),N} \right) \cong {\mathop{\rm Ext}\nolimits} _\mathcal{A}^n( {{\mathop{\rm Im}\nolimits} {d_1},{i^!}\left( N \right)} ) = 0,$$
thus ${\mathop{\rm Ext}\nolimits} _\mathcal{B}^n\left( {{C_0},N} \right) = 0$, that is, ${C_0} \in {}^ \bot {\rm{{\cal X}}}$. So we can obtain an exact sequence
$$0\longrightarrow P \stackrel{\varepsilon_0}\longrightarrow {{i_*}({T_0})}\oplus{{j_!}({Q_0})} \stackrel{\varepsilon_1}\longrightarrow {{i_*}({T_1})}\oplus {{j_!}({Q_1})} \stackrel{\varepsilon_2}\longrightarrow \dots,$$
by repeating the above process, where ${\mathop{\rm Im}\nolimits} {\kern 1pt} {\varepsilon _i} \in {}^ \bot {\rm{{\cal X}}}$ for all $i \geq 0$. Therefore, $\mathcal{X}$ is a Wakamatsu tilting subcategory of ${\rm{{\cal B}}}$.
\end{proof}

\begin{lemma}\rm\label{ll}
Let $f: \mathcal{A}\rightarrow \mathcal{B}$ and $g: \mathcal{B}\rightarrow \mathcal{A}$ be two functors between two abelian categories $\mathcal{A}$ and $\mathcal{B}$. Suppose that $(g, f)$ is a adjoint pair, both $f$ and $g$ are exact.
\begin{itemize}
  \item [(1)] If $\mathcal{X}$ is a self-orthogonal subcategory of $\mathcal{A}$ and $gf(\mathcal{X})\subseteq \mathcal{X}$, then $f(\mathcal{X})$ is a self-orthogonal subcategory of $\mathcal{B}$.
  \item [(2)] If $\mathcal{Y}$ is a self-orthogonal subcategory of $\mathcal{B}$ and $fg(\mathcal{Y})\subseteq \mathcal{Y}$, then $g(\mathcal{Y})$ is a self-orthogonal subcategory of $\mathcal{A}$.
\end{itemize}
\end{lemma}
\begin{proof}
(1) Let $f\left( {{M}} \right)$, $f\left( {{N}} \right) \in f\left( {\rm{{\cal X}}} \right)$. Then ${\mathop{\rm Ext}\nolimits} _{\mathcal{B}}^n( {f( {{M}} ),f( {{N}} )} ) \cong {\mathop{\rm Ext}\nolimits} _{\mathcal{A}}^n( {gf( {{M}} ),{N}})$. Since ${\rm{{\cal X}}}$ is a self-orthogonal subcategory of ${\rm{{\cal A}}}$ and $gf\left( {\rm{{\cal X}}} \right) \subseteq {\rm{{\cal X}}}$, so ${\mathop{\rm Ext}\nolimits} _{\mathcal{A}}^n( {gf( {{M}} ),{N}})=0.$ Hence $f(\mathcal{X})$ is a self-orthogonal subcategory of $\mathcal{B}$.

(2) It is similar to (1).
\end{proof}

\begin{lemma}\rm\label{rr}
Let $f: \mathcal{A}\rightarrow \mathcal{B}$ and $g: \mathcal{B}\rightarrow \mathcal{A}$ be two functors between two abelian categories $\mathcal{A}$ and $\mathcal{B}$. If $(g, f)$ is a adjoint pair, both $f$ and $g$ are exact, $g$ is fully faithful, $\mathcal{X}$ is a Wakamatsu tilting subcategory of $\mathcal{A}$ and $gf(^{\bot}\mathcal{X})\subseteq \mathcal{X}$, then $f(\mathcal{X})$ is a Wakamatsu tilting subcategory of $\mathcal{B}$.
\end{lemma}
\begin{proof}
Since $\mathcal{X}$ is a Wakamatsu tilting subcategory of $\mathcal{A}$, then $\mathcal{X}$ is a self-orthogonal subcategory of ${\rm{{\cal A}}}$. So $gf\left( {\rm{{\cal X}}} \right) \subseteq gf\left( {{}^ \bot {\rm{{\cal X}}}} \right) \subseteq {\rm{{\cal X}}}$. Hence $f\left( {\rm{{\cal X}}} \right)$ is a self-orthogonal subcategory of $\mathcal{B}$ by Lemma \ref{ll}. Next we show that ${X_{f\left( {\rm{{\cal X}}} \right)}}$ contains all projective objects of $\mathcal{B}$. Let $P$ be a projective object of $\mathcal{B}$, then $g\left( P \right)$ is a projective object of $\mathcal{A}$, hence $g\left( P \right) \in {X_{\rm{{\cal X}}}}$. Then there is an exact sequence
   $$\xymatrix{0\ar[r]&g(P)\ar[r]^{d_0}&{W_0}\ar[r]^{d_1}&{W_1}\ar[r] &\cdots}$$
with ${W_i} \in {\rm{{\cal X}}}$,  ${\mathop{\rm Im}\nolimits} {\kern 1pt} {d_i} \in {}^ \bot {\rm{{\cal X}}}$ for each $i \ge 0$. Since $f$ is exact, there is an exact sequence
$$\xymatrix{0\ar[r]&fg(P)\ar[r]^{f({d_0})}&f({W_0})\ar[r]^{f({d_1})}&f({W_1})\ar[r] &\cdots}.$$
Note that $fg\left( P \right) \cong P$, so we obtain an exact sequence
$$\xymatrix{0\ar[r]&P\ar[r]&f({W_0})\ar[r]^{f({d_1})}&f({W_1})\ar[r] &\cdots}.$$
Since ${\mathop{\rm Im}\nolimits} {\kern 1pt} f({d_i}) = {\mathop{\rm Ker}\nolimits} {\kern 1pt} f({d_{i + 1}})$, with f is exact, then $${\mathop{\rm Im}\nolimits} {\kern 1pt} f({d_i}) = {\mathop{\rm Ker}\nolimits} {\kern 1pt} f({d_{i + 1}}) \cong f{\mathop{\rm Ker}\nolimits} \left( {{d_{i + 1}}} \right).$$  Hence for each $f\left( M \right) \in f\left( {\rm{{\cal X}}} \right)$, we have $${\mathop{\rm Ext}\nolimits} _{\mathcal{B}}^n\left( {{\mathop{\rm Im}\nolimits} f\left( {{d_i}} \right), f\left( M \right)} \right)\cong {\mathop{\rm Ext}\nolimits} _{\mathcal{B}}^n\left( {f{\mathop{\rm Ker}\nolimits} \left( {{d_{i + 1}}} \right),f\left( M \right)} \right) \cong {\mathop{\rm Ext}\nolimits} _{\rm{{\cal A}}}^n\left( {gf{\mathop{\rm Ker}\nolimits} \left( {{d_{i + 1}}} \right),M} \right) = 0.$$ This show that $f\left( {\rm{{\cal X}}} \right)$ is a Wakamatsu tilting subcategory of $\mathcal{B}.$
\end{proof}

\begin{lemma}\rm\label{rrr}
Let $f: \mathcal{B}\rightarrow \mathcal{A}$, $g: \mathcal{A}\rightarrow \mathcal{B}$ and $h: \mathcal{B}\rightarrow \mathcal{A}$ be three functors between two abelian categories $\mathcal{A}$ and $\mathcal{B}$. If $(f, g, h)$ is a adjoint triple, $h$ is exact, $f$ is fully faithful, $\mathcal{X}$ is a Wakamatsu tilting subcategory of $\mathcal{A}$, $hg(\mathcal{X})\subseteq\mathcal{X}$, then $g(\mathcal{X})$ is a Wakamatsu tilting subcategory of $\mathcal{B}$.
\end{lemma}
\begin{proof}
Since $\mathcal{X}$ is a Wakamatsu tilting subcategory of $\mathcal{A}$, then $\mathcal{X}$ is a self-orthogonal subcategory of ${\rm{{\cal A}}}$. Hence $g\left( {\rm{{\cal X}}} \right)$ is a self-orthogonal subcategory of $\mathcal{B}$ by Lemma \ref{ll}. Next we show that ${X_{g\left( {\rm{{\cal X}}} \right)}}$ contains all projective objects of $\mathcal{B}$. Let $P$ be a projective object of $\mathcal{B}$, then $f\left( P \right)$ is a projective object of $\mathcal{A}$, hence $f\left( P \right) \in {X_{\rm{{\cal X}}}}$. Then there is an exact sequence
   $$\xymatrix{0\ar[r]&f(P)\ar[r]^{d_0}&{W_0}\ar[r]^{d_1}&{W_1}\ar[r] &\cdots}$$
with ${W_i} \in {\rm{{\cal X}}}$,  ${\mathop{\rm Im}\nolimits} {\kern 1pt} {d_i} \in {}^ \bot {\rm{{\cal X}}}$ for each $i \ge 0$. Since $g$ is exact, there is an exact sequence
$$\xymatrix{0\ar[r]&gf(P)\ar[r]^{g({d_0})}&g({W_0})\ar[r]^{g({d_1})}&g({W_1})\ar[r] &\cdots}.$$
Note that $P \cong gf\left( P \right)$, so we obtain an exact sequence
$$\xymatrix{0\ar[r]&P\ar[r]&g({W_0})\ar[r]^{g({d_1})}&g({W_1})\ar[r] &\cdots}.$$
Since ${\mathop{\rm Im}\nolimits} {\kern 1pt} g({d_i}) = {\mathop{\rm Ker}\nolimits} {\kern 1pt} g({d_{i + 1}})$, with g is exact, then $${\mathop{\rm Im}\nolimits} {\kern 1pt} g({d_i}) = {\mathop{\rm Ker}\nolimits} {\kern 1pt} g({d_{i + 1}}) \cong g{\mathop{\rm Ker}\nolimits} \left( {{d_{i + 1}}} \right).$$  Hence for each $g\left( M \right) \in g\left( {\rm{{\cal X}}} \right)$, we have $${\mathop{\rm Ext}\nolimits} _{\mathcal{B}}^n\left( {{\mathop{\rm Im}\nolimits} g\left( {{d_i}} \right), g\left( M \right)} \right)\cong {\mathop{\rm Ext}\nolimits} _{\mathcal{B}}^n\left( {g{\mathop{\rm Ker}\nolimits} \left( {{d_{i + 1}}} \right),g\left( M \right)} \right) \cong {\mathop{\rm Ext}\nolimits} _{\rm{{\cal A}}}^n\left( {{\mathop{\rm Ker}\nolimits} \left( {{d_{i + 1}}} \right),hg(M)} \right) = 0.$$ This show that $g\left( {\rm{{\cal X}}} \right)$ is a Wakamatsu tilting subcategory of $\mathcal{B}.$
\end{proof}

The following result shows that the converse of Theorem \ref{main1} holds true under certain
 conditions.
\begin{theorem}\label{main2}\rm
Let ($\mathcal{A}$, $\mathcal{B}$, $\mathcal{C}$) be a recollement of abelian categories as \rm{(\ref{recolle})}, $\mathcal{Y}$ is a Wakamatsu tilting subcategory of $\mathcal{B}$. If ${i^!}$ is exact, ${i_*}{i^!}\left( {{}^ \bot {\rm{{\cal Y}}}} \right) \subseteq {\rm{{\cal Y}}}$, ${j_*}{j^*}\left( {\rm{{\cal Y}}} \right) \subseteq {\rm{{\cal Y}}}$, then ${i^!}\left( {\rm{{\cal Y}}} \right)$ and ${j^*}\left( {\rm{{\cal Y}}} \right)$ are Wakamatsu tilting subcategories of $\mathcal{A}$ and $\mathcal{C}$, respectively.
\end{theorem}

\begin{proof}It follows from Lemma \ref{rr} and  \ref{rrr}.
\end{proof}

\section{Gluing weak support $\tau$-tilting subcategories }
For a subcategory $\mathcal{M}$ of an abelian category $\mathcal{A}$, the subcategory $\mathrm{Fac}(\mathcal{M})$ of $\mathcal{A}$ is defined as follow:

$\mathrm{Fac}(\mathcal{M}):=\{ C \in \mathcal{A} \mid $ there exists an epimorphism $M \rightarrow C \rightarrow 0 $, where $ M \in \mathcal{M}\}$.

Let $\mathcal{Y}$ be a subcategory of an abelian category $\mathcal{A}$ . A morphism $\varphi :A \to Y$, where $A$ is an object of $\mathcal{A}$ , is called
a left $\mathcal{Y}$-approximation of A, if $Y \in {\rm{{\cal Y}}}$ and for every $Y' \in {\rm{{\cal Y}}}$ , the induced sequence ${\mathop{\rm Hom}\nolimits}_\mathcal{A} \left( {Y, Y'} \right) \to {\mathop{\rm Hom}\nolimits}_\mathcal{A} \left( {A,Y'} \right) \to 0$
 of abelian groups is exact. We say that $\mathcal{Y}$ is a covariantly
finite subcategory of $\mathcal{A}$ if every object A of $\mathcal{A}$ admits a left $\mathcal{Y}$-approximation. Dually, the
notions of right $\mathcal{Y}$-approximations and contravariantly finite subcategories are defined.

\begin{definition}\label{support}{\rm \cite[Definition 3.1]{AST}}
Let $\mathcal{A}$ be an abelian category and $\mathcal{M}$ be a subcategory of $\mathcal{A}$. Then $\mathcal{M}$ is called a weak support $\tau$-tilting subcategory of $\mathcal{A}$ if it satisfies the following conditions:
 \begin{itemize}
  \item [(1)]  ${\mathop{\rm Ext}\nolimits} _{\rm{{\cal A}}}^1\left( {{\rm{{\cal M}}},{\mathop{\rm Fac}\nolimits} \left( {\rm{{\cal M}}} \right)} \right) = 0$.
  \item [(2)]   For any projective object $P$ in $\mathcal{A}$ , there exists an exact sequence $\xymatrix{P\ar[r]^{m}&M_{0} \ar[r]&{M_1}\ar[r]&0},$ such that ${M_0}, {M_1}\in {\rm{{\cal M}}}$ and $m$ is a left  $\mathcal{M}$-approximation of $P$.
\end{itemize}
\end{definition}
If furthermore $\mathcal{M}$ is a contravariantly finite subcategory of $\mathcal{A}$ , it is called a support $\tau$-tilting
subcategory of $\mathcal{A}$.



Our second main result is the following.
\begin{theorem}\label{main3}\rm
Let ($\mathcal{A}$, $\mathcal{B}$, $\mathcal{C}$) be a good recollement of abelian categories as \rm{(\ref{recolle})}, ${\rm{{\cal Z}}}'$ and ${\rm{{\cal Z}}}''$ are weak support $\tau$-tilting subcategories of $\mathcal{A}$ and $\mathcal{C}$.  Then $${\rm{{\cal Z}}} = \left\{ {Z \in {\rm{{\cal B}}}~|~{i^!}\left(Z \right) \in {\rm{{\cal Z}}}',{j^*}\left( Z \right) \in {\rm{{\cal Z}}}''} \right\}$$ is a weak support $\tau$-tilting subcategory of $\mathcal{B}$.
\begin{proof}
$\mathbf{Step ~1}$: Let $L$ $ \in {\mathop{\rm Fac}\nolimits} \left( {\rm{{\cal Z}}} \right)$, there is an exact sequence
$$\xymatrix{M_1\ar[r]&L\ar[r]&0}$$
with ${M_1}\in {\rm{{\cal Z}}}$. So we obtain two exact sequences
$$\xymatrix{i^!(M_1)\ar[r]&i^!(L)\ar[r]&0}$$ and
$$\xymatrix{j^*(M_1)\ar[r]&j^*(L)\ar[r]&0},$$
where ${i^!}\left( {{M_1}} \right) \in {\rm{{\cal Z}}}'$, ${j^*}\left( {{M_1}} \right)\in {\rm{{\cal Z}}}''.$ Hence ${i^!}\left( L \right) \in {\mathop{\rm Fac}\nolimits} \left( {{\rm{{\cal Z}}}'} \right),$ ${j^*}\left( L \right) \in {\mathop{\rm Fac}\nolimits} \left( {{\rm{{\cal Z}}}''} \right).$ Let $S \in {\rm{{\cal Z}}}.$ By Lemma \ref{7}, there are two exact sequences
$$\xymatrix{0\ar[r]&i_*i^!(S)\ar[r]&S\ar[r]&j_*j^*(S)\ar[r]&0}$$ and
$$\xymatrix{0\ar[r]&i_*i^!(L)\ar[r]&L\ar[r]&j_*j^*(L)\ar[r]&0}.$$
Hence we obtain the following two exact sequences
$$\xymatrix{\mathrm{Ext}_\mathcal{B}^1(j_*j^*(S),L)\ar[r]&\mathrm{Ext}_\mathcal{B}^1(S,L)\ar[r]&\mathrm{Ext}_\mathcal{B}^1(i_*i^!(S),L)}$$ and
$$\xymatrix{\mathrm{Ext}_\mathcal{B}^1(j_*j^*(S),i_*i^!(L))\ar[r]&\mathrm{Ext}_\mathcal{B}^1(j_*j^*(S),L)\ar[r]&\mathrm{Ext}_\mathcal{B}^1(j_*j^*(S),j_*j^*(L))}.$$
Then we have $${\mathop{\rm Ext}\nolimits} _{\rm{{\cal B}}}^1( {{i_*}{i^!}( S ),L} ) \cong {\mathop{\rm Ext}\nolimits} _{\rm{{\cal A}}}^1( {{i^!}( S ),{i^!}( L )}) = 0,$$
$${\mathop{\rm Ext}\nolimits} _{\rm{{\cal B}}}^1( {{j_*}{j^*}( S ),{i_*}{i^!}( L )}) \cong {\mathop{\rm Ext}\nolimits} _{\rm{{\cal A}}}^1( {{i^*}{j_*}{j^*}( S ),{i^!}( L )}) = 0$$ and $${\mathop{\rm Ext}\nolimits} _{\rm{{\cal B}}}^1\left( {{j_*}{j^*}\left( S \right),{j_*}{j^*}\left( L \right)} \right) \cong {\mathop{\rm Ext}\nolimits} _{\rm{{\cal C}}}^1\left( {{j^*}\left( S \right),{j^*}{j_*}{j^*}\left( L \right)} \right) \cong {\mathop{\rm Ext}\nolimits} _{\rm{{\cal C}}}^1\left( {{j^*}\left( S \right),{j^*}\left( L \right)} \right) = 0,$$ by Lemma \ref{7}, so ${\mathop{\rm Ext}\nolimits} _{\rm{{\cal B}}}^1\left( {{j_*}{j^*}\left( S \right),L} \right) = 0,$ hence ${\mathop{\rm Ext}\nolimits} _{\rm{{\cal B}}}^1\left( {S,L} \right) = 0,$ then ${\mathop{\rm Ext}\nolimits} _{\rm{{\cal B}}}^1\left( {{\rm{{\cal Z}}},{\mathop{\rm Fac}\nolimits} \left( {\rm{{\cal Z}}} \right)} \right) = 0$.

$\mathbf{Step ~2}$: Let $P$ be a projective object of ${\rm{{\cal B}}}$, then ${i^*}\left( P \right)$ and ${j^*}\left( P \right)$ are projective objects of ${\rm{{\cal A}}}$ and ${\rm{{\cal C}}}$, respectively, by Lemma \ref{7}. Then there are two exact sequences
$$\xymatrix{i^*(P)\ar[r]^{m}&Z^0\ar[r]^{d_0}&Z^1\ar[r]&0}$$ and
$$\xymatrix{j^*(P)\ar[r]^{n}&Y^0\ar[r]^{q_0}&Y^1\ar[r]&0},$$
with ${Z^0}$, ${Z^1} \in {\rm{{\cal Z}}}'$, ${Y^0}$, ${Y^1} \in {\rm{{\cal Z}}}''$ and $m$ is a left $\mathcal{Z'}$-approximation of ${i^*}\left( P \right)$, $n$ is a left $\mathcal{Z''}$-approximation of ${j^*}\left( P \right)$. So there are two exact sequences
$$\xymatrix{{i_*}{i^*}(P)\ar[r]^{i_*(m)}&i_*(Z^0)\ar[r]^{i_*({d_0})}&{i_*(Z^1)}\ar[r]&0}$$ and
$$\xymatrix{{j_!}{j^*}(P)\ar[r]^{j_!(n)}&j_!(Y^0)\ar[r]^{j_!({q_0})}&{j_!(Y^1)}\ar[r]&0}.$$
Note that there exists a split exact sequence
\begin{equation}\label{4.2}\xymatrix{0\ar[r]&j_!j^*(P)\ar[r]^{h}&P\ar[r]^{g}&i_*i^*(P)\ar[r]&0}\end{equation}
by Lemma \ref{7} and Remark \ref{Yu}. That means there is a isomorphism $f:P \to {j_!}{j^*}\left( P \right) \oplus {i_*}{i^*}\left( P \right)$. Then there is an exact sequence
$$ P \xrightarrow[]{\Big( {\begin{smallmatrix}
{{j_!}(n)}&{0}\\
{0}&{{i_*}(m)}
\end{smallmatrix}}\Big)f} {{j_!}({Y^0})}\oplus{{i_*}({Z^0})} \xrightarrow[]{\Big( {\begin{smallmatrix}
{{j_!}({q_0})}&{0}\\
{0}&{{i_*}({d_0})}
\end{smallmatrix}}\ \Big)} {{j_!}({Y^1})}\oplus {{i_*}({Z^1})}.$$
Since ${i^!}\left( {{j_!}\left( {{Y^t}} \right) \oplus {i_*}\left( {{Z^t}} \right)} \right) \cong {Z^t},$ ${j^*}\left( {{j_!}\left( {{Y^t}} \right) \oplus {i_*}\left( {{Z^t}} \right)} \right) \cong {Y^t},$ so ${j_!}\left( {{Y^t}} \right) \oplus {i_*}\left( {{Z^t}} \right) \in {\rm{{\cal Z}}},~t=0,1.$ Next, we show that $\Big( {\begin{smallmatrix}
{{j_!}(n)}&{0}\\
{0}&{{i_*}(m)}
\end{smallmatrix}} \Big)f$
is a left $\mathcal{Z}$-approximation of $P$. We need to show that
$$ \mathcal{B}(j_!(Y^0)\oplus i_*(Z^0),Z) \xrightarrow[]{\mathcal{B}\Big(\Big( {\begin{smallmatrix}
{{j_!}(n)}&{0}\\
{0}&{{i_*}(m)}
\end{smallmatrix}}\Big)f, ~Z\Big)} \mathcal{B}(P,Z) \xrightarrow[]{}0$$
is exact for any $Z \in {\rm{{\cal Z}}}$. Since $${\mathcal{B}\left( {\left( {\begin{smallmatrix}
{{j_!}(n)}&{0}\\
{0}&{{i_*}(m)}
\end{smallmatrix}} \right)f} , Z\right)} = \mathcal{B}(f, Z)\circ\mathcal{B}\Big({\left( {\begin{smallmatrix}
{{j_!}(n)}&{0}\\
{0}&{{i_*}(m)}
\end{smallmatrix}} \right), ~Z\Big)},$$
it is enough to show that $${\mathcal{B}\left(\left( {\begin{smallmatrix}
{{j_!}(n)}&{0}\\
{0}&{{i_*}(m)}
\end{smallmatrix}} \right), Z\right)}:{{{\rm{{\cal B}}}}\left( {{j_!}\left( {{Y^0}} \right) \oplus {i_*}\left( {{Z^0}} \right),Z} \right) \to {\rm{{\cal B}}}}\left( {{j_!}{j^*}\left( P \right) \oplus {i_*}{i^*}\left( P \right),Z} \right)$$ is epimorphism. And there is a commutative diagram
\begin{center}
$\xymatrix{
\mathcal{B}( {j_!}( {{Y^0}} ) \oplus {i_*}( {{Z^0}} ), Z )\ar[r]\ar[d]^{\mathcal{B}\Big(\Big( {\begin{smallmatrix}
{{j_!}(n)}&{0}\\
{0}&{{i_*}(m)}
\end{smallmatrix}}\Big), Z\Big)} &\mathcal{B}(j_!(Y^0), Z)\oplus \mathcal{B}(i_*(Z^0),Z)\ar[d]^{\Big( {\begin{smallmatrix}
{\mathcal{B}({j_!}(n), Z)}&{0}\\
{0}&{\mathcal{B}({i_*}(m), Z)}
\end{smallmatrix}}\Big)}\\
\mathcal{B}(j_!j^*(P)\oplus i_*i^*(P),Z)\ar[r]&\mathcal{B}( {j_!}{j^*}( P ),Z ) \oplus \mathcal{B}( {i_*}{i^*}( P ),Z)} $
\end{center}
it is enough to show that $\left( {\begin{smallmatrix}
\mathcal{B}({{j_!}(n),Z)}&{0}\\
{0}&{\mathcal{B}({i_*}(m),Z)}
\end{smallmatrix}} \right)$
is epimorphism. Since $\left( {{i_*},{i^!}} \right)$ and $\left( {{j_!},{j^*}} \right)$ are adjoint pairs, there is a commutative diagram
\begin{center}
$\xymatrix{
\mathcal{C}(  Y^0 ,j^*(Z) )\oplus  \mathcal{A}( Z^0 ,i^!(Z) )\ar[r]^{\cong}\ar[d]^{\Big({\begin{smallmatrix}
{\mathcal{C}(n, j^*(Z))}&{0}\\
{0}&{\mathcal{A}(m, i^!(Z))}
\end{smallmatrix}}\Big)} &\mathcal{B}(j_!(Y^0), Z)\oplus \mathcal{B}(i_*(Z^0), Z)\ar[d]^{\Big( {\begin{smallmatrix}
{\mathcal{B}({j_!}(n), Z)}&{0}\\
{0}&{\mathcal{B}({i_*}(m),Z)}
\end{smallmatrix}}\Big)}\\
\mathcal{C}(j^*(P),j^*(Z))\oplus \mathcal{A}(i^*(P),i^!(Z)) \ar[r]^{\cong}&\mathcal{B}( {j_!}{j^*}( P ), Z ) \oplus \mathcal{B}( {i_*}{i^*}( P ), Z)} $
\end{center}
It is enough to show that $\left( {\begin{smallmatrix}
{{{\mathcal{C}(n, j^*( Z )) }}}&{0}\\
{0}&{{{\mathcal{A}(m,i^!( Z )) }}}
\end{smallmatrix}} \right)$
is epimorphism. Since $m$ is a left $\mathcal{Z'}$-approximation of ${i^*}\left( P \right)$ and $n$ is a left $\mathcal{Z''}$-approximation of ${j^*}\left( P \right)$, so
$$\mathcal{A}( Z^0,i^!( Z )) \xrightarrow[]{{\mathcal{A}(m, i^!(Z))}} \mathcal{A}( {i^*}( P ),i^!( Z ) ) \xrightarrow[]{} 0,$$
$$\mathcal{C}( Y^0,j^*( Z )) \xrightarrow[]{{\mathcal{C}(n, j^*(Z))}} \mathcal{C}( {j^*}( P ),j^*( Z ) ) \xrightarrow[]{} 0$$
are exact with ${i^!}\left( Z \right) \in {\rm{{\cal Z}}}'$ and ${j^*}\left( Z \right) \in {\rm{{\cal Z}}}''$, then $\left( {\begin{smallmatrix}
{{{\mathcal{C}(n, j^*( Z )) }}}&{0}\\
{0}&{{{\mathcal{A}(m,i^!( Z )) }}}
\end{smallmatrix}} \right)$
is an epimorphism.
\end{proof}
\end{theorem}
\begin{lemma}\rm \label{vv}
Let $f: \mathcal{B}\rightarrow \mathcal{A}$, $g: \mathcal{A}\rightarrow \mathcal{B}$ and $h: \mathcal{B}\rightarrow \mathcal{A}$ be three functors between two abelian categories $\mathcal{A}$ and $\mathcal{B}$. Suppose that $(f, g, h)$ is a adjoint triple, all $f$, $g$ and $h$ are exact.
 \begin{itemize}
\item [(1)]If $g$ is fully faithful, $\mathcal{Y}$ is a weak support $\tau$-tilting subcategory of $\mathcal{B}$ and $gf({\mathcal{Y}})\subseteq \mathcal{Y}$, then $f(\mathcal{Y})$ is a weak support $\tau$-tilting subcategory of $\mathcal{A}$.
\item [(2)] If $f$ is fully faithful, $\mathcal{Y}$ is a weak support $\tau$-tilting subcategory of $\mathcal{A}$ and $hg({\mathcal{Y}})\subseteq \mathcal{Y}$, then $g(\mathcal{Y})$ is a weak support $\tau$-tilting subcategory of $\mathcal{B}$.
 \end{itemize}

\end{lemma}
\begin{proof}
(1) $\mathbf{Step ~1}$: Let $K \in {\mathop{\rm Fac}\nolimits} \left( {f\left( {\rm{{\cal Y}}} \right)} \right)$, there is an exact sequence
$$\xymatrix{f({M_1})\ar[r]&K\ar[r] &0}$$
where ${{M_1}} $ $ \in {\rm{{\cal Y}}} $. Since $g$ is exact, there is an exact sequence
$$\xymatrix{gf({M_1})\ar[r]&g(K)\ar[r] &0}.$$ Since $gf\left( {\rm{{\cal Y}}} \right) \subseteq {\rm{{\cal Y}}}$, then $g\left( K \right)$ $ \in {\mathop{\rm Fac}\nolimits} \left( {\rm{{\cal Y}}} \right)$. Let ${{M_2}}$ $ \in {\rm{{\cal Y}}} $, then we have
$${\mathop{\rm Ext}\nolimits} _{\rm{{\cal A}}}^1\left( {f\left( {{M_2}} \right),K} \right) \cong {\mathop{\rm Ext}\nolimits} _{\rm{{\cal A}}}^1\left( {{M_2},g\left( K \right)} \right) = 0.$$
That is, ${\mathop{\rm Ext}\nolimits} _{\rm{{\cal A}}}^1\left( {{f\left( {\rm{{\cal Y}}} \right),~{\mathop{\rm Fac}\nolimits} \left( {f\left( {\rm{{\cal Y}}} \right)} \right)}} \right) = 0.$

$\mathbf{Step ~2}$: Let $P$ be a projective object of $\mathcal{A}$, then $g\left( P \right)$ is a projective object of ${\rm{{\cal B}}}$. Hence there is an exact sequence
$$\xymatrix{g(P)\ar[r]^m&X^0\ar[r] &X^1\ar[r]&0},$$ where ${X^0}, {X^1}\in \mathcal{Y}$, $m$ is a left $\mathcal{Y}$-approximation of $g\left( P \right)$. Then there is an exact sequence
$$\xymatrix{fg(P)\ar[r]^{f(m)}&f(X^0)\ar[r]&{f(X^1)}\ar[r]&0}.$$
Since $P\overset{t}\cong fg(P)$, there is an exact sequence
$$P\overset{f(m)t}\longrightarrow f({X^0}) \longrightarrow {f({X^1})} \longrightarrow 0.$$
Let ${ Y  \in {\rm{{\cal Y}}} }$, since $\left( f,g \right)$ is an adjoint pair, there is a commutative diagram
\begin{center}
$\xymatrix{
\mathcal{A}(f(X^0),f(Y))\ar[r]\ar[d]^{\cong}&\mathcal{A}(fg(P), f(Y))\ar[r]^{\cong}\ar[d]^{\cong}&\mathcal{A}(P, f(Y))\\
\mathcal{B}(X^0,gf(Y))\ar[r]&\mathcal{B}(g(P),gf(Y))}$
\end{center}
Since ${gf\left( {\rm{{\cal Y}}} \right) \subseteq {\rm{{\cal Y}}}}$, there is an exact sequence
$$\xymatrix{\mathcal{B}(X^0,gf(Y))\ar[r]&\mathcal{B}(g(P),gf(Y))\ar[r]&0},$$
hence
$$\xymatrix{\mathcal{A}(f(X^0), f(Y))\ar[r]&\mathcal{A}(P, f(Y))\ar[r]&0}$$
is exact. Then $f(m)t$ is a left $f(\mathcal{Y})$-approximation of $P$.
Therefore, $f\left( {\rm{{\cal Y}}} \right)$ is a weak support $\tau$-tilting subcategory of ${\rm{{\cal A}}}$.

(2)  It is similar to (1).

\end{proof}

The following result shows that the converse of Theorem \ref{main3} holds true under certain
 conditions.
\begin{theorem}\label{main4}\rm
Let ($\mathcal{A}$, $\mathcal{B}$, $\mathcal{C}$) be a good recollement of abelian categories as \rm{(\ref{recolle})}, ${\rm{{\cal Y}}}$ is a weak support $\tau$-tilting subcategory of ${\rm{{\cal B}}}$. If ${i_*}{i^*}\left( {\rm{{\cal Y}}} \right) \subseteq {\rm{{\cal Y}}}$, ${j_*}{j^*}\left( {\rm{{\cal Y}}} \right) \subseteq {\rm{{\cal Y}}}$, then ${i^*}\left( {\rm{{\cal Y}}} \right)$ and ${j^*}\left( {\rm{{\cal Y}}} \right)$ are weak support $\tau$-tilting subcategories of ${\rm{{\cal A}}}$ and ${\rm{{\cal C}}}, $ respectively.
\begin{proof} It follows from Lemma \ref{vv}.
\end{proof}
\end{theorem}

\section{Application}

Asadollahi, Sadeghi and Treffinger \cite{AST} established a bijection between $\tau\textrm{-cotorsion torsion triples}$ and $\textrm{support}~\tau\textrm{-tilting subcategories}$. In this section, we consider the relationship of $\tau$-cotorsion torsion triples in  ($\mathcal{A}$, $\mathcal{B}$, $\mathcal{C}$).
\begin{lemma}\label{main5}\rm
Let ($\mathcal{A}$, $\mathcal{B}$, $\mathcal{C}$) be a good recollement of abelian categories as \rm{(\ref{recolle})}, ${\rm{{\cal Z}}}'$ and ${\rm{{\cal Z}}}''$ are contravariantly finite subcategories of $\mathcal{A}$ and $\mathcal{C}$. Then $${\rm{{\cal Z}}} = \left\{ {Z \in {\rm{{\cal B}}}~|~{i^*}\left( Z \right) \in {\rm{{\cal Z}}}',{j^*}\left( Z \right) \in {\rm{{\cal Z}}}''} \right\}$$ is contravariantly finite subcategory of $\mathcal{B}$.
\begin{proof}
Let $B\in \mathcal{B}$, then $i^{!}(B)\in \mathcal{A}$ and $j^*(B)\in \mathcal{C}$. So there are two morphisms $\varphi: W'\rightarrow i^!(B)$ and $\psi: R''\rightarrow j^*(B),$ where $W'\in \mathcal{Z}',$ $R''\in \mathcal{Z}''$, $\varphi$ is a right $\mathcal{Z}'$-approximation of $i^!(B)$ and $\psi$ is a right $\mathcal{Z}''$-approximation of $j^*(B)$. So there are two morphisms $i_*(\varphi): i_*(W')\rightarrow i_*i^!(B),$ $j_*(\psi):j_*(R'')\rightarrow j_*j^*(B).$
Note that there is a split exact sequence \begin{equation}\label{5.3}0\rightarrow i_*i^!(B)\rightarrow B\rightarrow j_*j^*(B)\rightarrow 0\end{equation}  by Lemma \ref{7} and Remark \ref{Yu}. So there is an isomorphism $l: j_*j^*(B)\oplus i_*i^!(B)\rightarrow B$. Therefore there is a morphism $ l\big(\begin{smallmatrix} &{j_*(\psi)} &{0}\\ &{0} &{i_*(\varphi)}\end{smallmatrix}\big): j_*(R'')\oplus i_*(W')\rightarrow B.$  Since $i^!(j_*(R'')\oplus i_*(W'))\cong W',$ $j^*(j_*(R'')\oplus i_*(W'))\cong R'',$ so $j_*(R'')\oplus i_*(W')\in \mathcal{Z}.$ Fix $F\in \mathcal{Z}$. Next we show that the sequence $$\mathcal{B}(F, j_*(R'')\oplus i_*(W'))\rightarrow \mathcal{B}(F,B)\rightarrow 0$$ is exact or $\mathcal{B}\Big(F,l\Big(\begin{smallmatrix}&{j_*(\psi)}&{0}\\&{0}&{i_*(\varphi)} \end{smallmatrix}\Big)\Big)$ is surjective. Since $$\mathcal{B}\Big(F, l\Big(\begin{smallmatrix}&{j_*(\psi)}&{0}\\&{0}&{i_*(\varphi)} \end{smallmatrix}\Big)\Big)=\mathcal{B}(F,l)\circ\mathcal{B}\Big(F,\Big(\begin{smallmatrix}&{j_*(\psi)}&{0}\\&{0}&{i_*(\varphi)} \end{smallmatrix}\Big)\Big),$$ it is enough to show that $\mathcal{B}\Big(F, \Big(\begin{smallmatrix}&{j_*(\psi)}&{0}\\&{0}&{i_*(\varphi)} \end{smallmatrix}\Big)\Big)$ is surjective. There is a commutative diagram
\begin{center}
$\xymatrix{
\mathcal{B}(F, j_*(R'')\oplus i_*(W'))\ar[r]^{\cong}\ar[d]^{\mathcal{B}\Big(F,\Big(\begin{smallmatrix}&{j_*(\psi)}&{0}\\&{0}&{i_*(\varphi)} \end{smallmatrix}\Big)\Big)}&\mathcal{B}( F, {j_*}( R'' ) ) \oplus\mathcal{B}( F, {i_*}( W' ))\ar[d]^{\Big(\begin{smallmatrix}&{\mathcal{B}(F,j_*(\psi))}&{0}\\&{0}&{\mathcal{B}(F,i_*(\varphi))} \end{smallmatrix}\Big)}\\
\mathcal{B}( F, {j_*}{j^*}(B) \oplus {i_*}{i^!}( B ) )\ar[r]^{\cong} &\mathcal{B}(F, j_*j^*(B))\oplus \mathcal{B}(F, i_*i^!(B))} $
\end{center}
so we just need to show that $\Big(\begin{smallmatrix}&{\mathcal{B}(F, j_*(\psi))}&{0}\\&{0}&{\mathcal{B}(F, i_*(\varphi))} \end{smallmatrix}\Big)$ is surjective.
Since there is a commutative diagram
\begin{center}
$\xymatrix{
\mathcal{C}( j^*(F), R'' ) \oplus\mathcal{A}( i^*(F), W')\ar[r]\ar[d]^{\cong} &\mathcal{C}(j^*(F), j^*(B))\oplus \mathcal{A}(i^*(F), i^!(B))\ar[d]^{\cong} \\
 \mathcal{B}( F, {j_*}( R'' ) ) \oplus \mathcal{B}( F, {i_*}( W' ))\ar[r]&\mathcal{B}(F, j_*j^*(B))\oplus \mathcal{B}(F, i_*i^!(B))}$
\end{center}
with $j^*(F)\in \mathcal{Z}''$, $i^*(F)\in \mathcal{Z}'$, $\varphi$ is a right $Z'$-approximation of $i^!(B)$ and $\psi$ is a right $Z''$-approximation of $j^*(B).$ So $\Big(\begin{smallmatrix}&{\mathcal{B}(F, j_*(\psi))}&{0}\\&{0}&{\mathcal{B}(F, i_*(\varphi))} \end{smallmatrix}\Big)$ is surjective.
\end{proof}
\end{lemma}
By Theorem \ref{main3} and Lemma \ref{main5}, we have the following.
\begin{corollary}\label{xx}
\rm Let ($\mathcal{A}$, $\mathcal{B}$, $\mathcal{C}$) be a good recollement of abelian categories as \rm{(\ref{recolle})}, ${\rm{{\cal Z}}}'$ and ${\rm{{\cal Z}}}''$ are support $\tau$-tilting subcategories of $\mathcal{A}$ and $\mathcal{C}$. Then $${\rm{{\cal Z}}} = \left\{ {Z \in {\rm{{\cal B}}}~|~{i^*}\left( Z \right) \in {\rm{{\cal Z}}}', i^!(Z)\in Z', {j^*}\left( Z \right) \in {\rm{{\cal Z}}}''} \right\}$$ is a support $\tau$-tilting subcategory of $\mathcal{B}$.
\end{corollary}
\begin{lemma}\rm\label{000}
Let $f: \mathcal{B}\rightarrow \mathcal{A}$, $g: \mathcal{A}\rightarrow \mathcal{B}$ and $h: \mathcal{B}\rightarrow \mathcal{A}$ be three functors between two abelian categories $\mathcal{A}$ and $\mathcal{B}$. Suppose that $(f, g, h)$ is a adjoint triple.
\begin{itemize}
\item [(1)]If $g$ is fully faithful, $\mathcal{Y}$ is a contravariantly finite subcategory of $\mathcal{B}$ and $gh({\mathcal{Y}})\subseteq \mathcal{Y}$, $gf({\mathcal{Y}})\subseteq \mathcal{Y}$, then $f(\mathcal{Y})$ is a contravariantly finite subcategory of $\mathcal{A}$.
\item [(2)]If $h$ is fully faithful, $\mathcal{Y}$ is a contravariantly finite subcategory of $\mathcal{A}$ and $fg({\mathcal{Y}})\subseteq \mathcal{Y}$, then $g(\mathcal{Y})$ is a contravariantly finite subcategory of $\mathcal{B}$.
\end{itemize}
\end{lemma}
\begin{proof}
(1) Let $A\in \mathcal{A}$, then $g(A)\in \mathcal{B}$. Since $\mathcal{Y}$ is a contravariantly finite subcategory of ${\rm{{\cal B}}}$, so there is a morphism $t: Z\rightarrow g(A)$ with $Z\in \mathcal{Y}$ and $t$ is right $\mathcal{Y}$-approximation of $g(A)$. Hence there is a morphism $h(t): h(Z)\rightarrow hg(A)$. Since $hg(A)\overset{k}\cong A $, $h(Z)\cong fgh(Z)$ and $gh(\mathcal{Y})\subseteq \mathcal{Y}$, $h(Z)\in f(\mathcal{Y})$ with $f(\mathcal{Y})$ is closed to isomorphism. Let $f(Y)\in f(\mathcal{Y})$, there is a commutative diagram
\begin{center}
$\xymatrix{
 \mathcal{A}( f(Y), h(Z) )\ar[r]\ar[d]^{\cong} &\mathcal{A}(f(Y), hg(A))  \ar[d]^{\cong} \xrightarrow[]{\cong}\mathcal{A}(f(Y), A)\\
\mathcal{B}( gf(Y), Z )\ar[r]&\mathcal{B}(gf(Y), g(A))}$
\end{center}
So we know that $$\mathcal{A}(f(Y), h(Z))\rightarrow\mathcal{A}(f(Y), A)\rightarrow 0$$ is exact since $$\mathcal{B}(gf(Y), Z)\rightarrow\mathcal{B}(gf(Y), g(A))\rightarrow 0$$ is exact, with $gf(Y)\in \mathcal{Y}$ and $t$ is a right $\mathcal{Y}$-approximation of $g(A)$, hence $kh(t)$ is a right $f(\mathcal{Y})$-approximation of $A$. This shows that $f\left( {\rm{{\cal Y}}} \right)$ is a contravariantly finite subcategory of ${\rm{{\cal A}}}.$
(2)  It is similar to (1).

\end{proof}

By Theorem \ref{main4} and Lemma \ref{000}, we have the following
\begin{corollary} \label{bb}
\rm Let ($\mathcal{A}$, $\mathcal{B}$, $\mathcal{C}$) be a good recollement of abelian categories as \rm{(\ref{recolle})}, ${\rm{{\cal Y}}}$ is a support $\tau$-tilting subcategory of ${\rm{{\cal B}}}$. If  ${i_*}{i^!}\left( {\rm{{\cal Y}}} \right) \subseteq {\rm{{\cal Y}}}$, $i_*i^*(\mathcal{Y})\subseteq \mathcal{Y},$  ${j_!}{j^*}\left( {\rm{{\cal Y}}} \right) \subseteq {\rm{{\cal Y}}}$, ${j_*}{j^*}\left( {\rm{{\cal Y}}} \right) \subseteq {\rm{{\cal Y}}}$, then ${i^*}\left( {\rm{{\cal Y}}} \right)$ and ${j^*}\left( {\rm{{\cal Y}}} \right)$ are support $\tau$-tilting subcategories of ${\rm{{\cal A}}}$ and ${\rm{{\cal C}}}, $ respectively.
\end{corollary}
Let $\mathcal{H}$ be a subcategory of an abelian category $\mathcal{A}$ and $n$ be a non-negative integer.

Set
$$\mathcal{H}^{\perp_n}:=\{A\in \mathcal{A}\mid \mathrm{Ext}_{\mathcal{A}}^n(\mathcal{H},A)=0\},$$
$$^{\perp_n}\mathcal{H}:=\{A\in \mathcal{A}\mid \mathrm{Ext}_{\mathcal{A}}^n(A, \mathcal{H})=0\}.$$
Note that $\mathrm{Ext}^0$ is just the usual $\mathrm{Hom}$-functor. Let us recall the concept of a $\tau$-cotorsion torsion triple of abelian categories from \cite {AST}
\begin{definition}\label{support}{\rm \cite[Definition 1.2]{AST}}
Let $\mathcal{A}$ be an abelian category. A triple of full subcategories$(\mathcal{L},\mathcal{D},\mathcal{F})$ of $\mathcal{A}$ is called a $\tau$-cotorsion torsion triple if
 \begin{itemize}
  \item [(1)] $\mathcal{L}={}^{\perp_1}\mathcal{D}.$
  \item [(2)] For every projective object $P\in \mathcal{A}$, there exists an exact sequence $$P\xrightarrow[]{f} D\xrightarrow[]{} C\xrightarrow[]{} 0,$$ where $D\in \mathcal{L}\cap \mathcal{D},~ C\in \mathcal{L}$ and $f$ is a left $\mathcal{D}$-approximation.
  \item [(3)] $\mathcal{L}\cap \mathcal{D}$ is a contravariantly finite subcategory of $\mathcal{A}$.
  \item [(4)] $(\mathcal{D},\mathcal{F})$ is a torsion pair in $\mathcal{A}$.
\end{itemize}
\end{definition}
Asadollahi, Sadeghi and Treffinger \cite{AST} established a bijections between $\tau\textrm{-cotorsion torsion triples}$ and $\textrm{support}~\tau\textrm{-tilting subcategories}$.

\begin{lemma}\label{uu}\rm\cite[Theorem 5.7]{AST}
Let $\mathcal{A}$ be an abelian category. Then there are bijections
$$\Phi: \{\textrm{support}~\tau\textrm{-tilting subcategories}\}\rightarrow \{\tau\textrm{-cotorsion torsion triples}\}$$
$$\mathcal{L}\mapsto ({}^{\perp_1}(\mathrm{Fac}(\mathcal{L})),\mathrm{Fac}(\mathcal{L}),\mathcal{L}^{\perp_0})$$
$$\Psi: \{\tau\textrm{-cotorsion torsion triples}\}\rightarrow \{\textrm{support}~\tau\textrm{-tilting subcategories}\}$$
$$(\mathcal{L}, \mathcal{D}, \mathcal{F})\mapsto \mathcal{L}\cap \mathcal{D}$$
which are mutually inverse.
\end{lemma}

\begin{proposition}\label{p1}
\rm Let ($\mathcal{A}$, $\mathcal{B}$, $\mathcal{C}$) be a good recollement of abelian categories as \rm{(\ref{recolle})}, $(\mathcal{E},\mathcal{F},\mathcal{G})$ and $(\mathcal{E}',\mathcal{F}',\mathcal{G}')$ are $\tau\textrm{-cotorsion torsion triples}$ of $\mathcal{A}$ and $\mathcal{C}$, respectively. Define
$$\mathcal{Z}=\{Z\in \mathcal{B}~|~i^*(Z)\in \mathcal{E}\cap \mathcal{F}, i^!(Z)\in \mathcal{E}\cap \mathcal{F}, j^*(Z)\in \mathcal{E}'\cap \mathcal{F}'\}$$
Then $(^{\perp_1}(\mathrm{Fac}(\mathcal{Z})), \mathrm{Fac}(\mathcal{Z}), \mathcal{Z}^{\perp_0})$ is a $\tau$-cotorsion torsion triple of $\mathcal{B}.$
\end{proposition}
\proof
By Lemma \ref{uu}, we know that $\mathcal{E}\cap \mathcal{F}$ and $\mathcal{E}'\cap \mathcal{F}'$ are $\textrm{support}~\tau\textrm{-tilting subcategories}$ of $\mathcal{A}$ and $\mathcal{C}$, respectively. Then $\mathcal{Z}$ is a $\textrm{support}~\tau\textrm{-tilting subcategory}$ of $\mathcal{B}$ by Corollary \ref{xx}. And hence $(^{\perp_1}(\mathrm{Fac}(\mathcal{Z})), \mathrm{Fac}(\mathcal{Z}), \mathcal{Z}^{\perp_0})$ is a $\tau$-cotorsion torsion triple of $\mathcal{B}$ by Lemma \ref{uu}.
\qed
\begin{proposition}\label{p2}
\rm Let ($\mathcal{A}$, $\mathcal{B}$, $\mathcal{C}$) be a good recollement of abelian categories as \rm{(\ref{recolle})}, $(\mathcal{E},\mathcal{F},\mathcal{G})$ is a $\tau$-cotorsion torsion triple of $\mathcal{B}$. If $i_*i^!(\mathcal{E}\cap \mathcal{F})\subseteq \mathcal{E}\cap \mathcal{F}$, $i_*i^*(\mathcal{E}\cap \mathcal{F})\subseteq \mathcal{E}\cap \mathcal{F}$, $j_!j^*(\mathcal{E}\cap \mathcal{F})\subseteq \mathcal{E}\cap \mathcal{F}$, $j_*j^*(\mathcal{E}\cap \mathcal{F})\subseteq \mathcal{E}\cap \mathcal{F}$, then
$$\clubsuit=(^{\perp_1}(\mathrm{Fac}(i^*(\mathcal{E}\cap \mathcal{F}))), \mathrm{Fac}(i^*(\mathcal{E}\cap \mathcal{F})), (i^*(\mathcal{E}\cap \mathcal{F}))^{\perp_0})$$ and $$\spadesuit=(^{\perp_1}(\mathrm{Fac}(j^*(\mathcal{E}\cap \mathcal{F}))), \mathrm{Fac}(j^*(\mathcal{E}\cap \mathcal{F})), (j^*(\mathcal{E}\cap \mathcal{F}))^{\perp_0})$$  are $\tau$-cotorsion torsion triples of $\mathcal{A}$ and $\mathcal{C}$, respectively.
\end{proposition}

\proof
By Lemma \ref{uu}, we know that $\mathcal{E}\cap \mathcal{F}$ is a $\textrm{support}~\tau\textrm{-tilting subcategory}$ of $\mathcal{B}$. Then $i^*(\mathcal{E}\cap \mathcal{F})$ and $j^*(\mathcal{E}\cap \mathcal{F})$ are support~$\tau$-tilting subcategories of $\mathcal{A}$ and $\mathcal{C}$ by Corollary \ref{bb}, respectively. And hence $\clubsuit$ and $\spadesuit$ are $\tau$-cotorsion torsion triples of $\mathcal{A}$ and $\mathcal{C}$ by Lemma \ref{uu}, respectively.
\qed

\section{Examples} \label{sect:examp}

\def\kk{\Bbbk} \def\Q{\mathcal{Q}} \def\I{\mathcal{I}}
\def\M{\mathfrak{M}} \def\X{\mathfrak{X}} \def\bfS{\mathbf{S}}
\def\rmark{{\color{red}\pmb{\circ}}} \def\bmark{{\color{blue}\bullet}}
\def\PC{\mathrm{PC}} \def\CC{\mathrm{CC}} \def\Str{\mathrm{Str}} \def\Band{\mathrm{Band}}
\def\ind{\mathrm{ind}}
\def\source{\mathfrak{s}} \def\target{\mathfrak{t}}

In this section, we provide some examples for Theorems \ref{main1}, \ref{main2}, \ref{main3}, and \ref{main4}.
Let $A=\kk\Q_A/\I_A$ be a finite-dimensional algebra over an algebraically closed field $\kk$ given by the following bound quiver $(\Q_A,\I_A)$:
\[ \Q_A = \ \
\xymatrix{
  1 \ar[r]^{\alpha}
& 2 \ar@/^0.5pc/[r]^{\beta_1}
    \ar@/_0.5pc/[r]_{\beta_2}
& 3 \ar[r]^{\gamma}
& 4,
}
\ \
\I_A = \langle \alpha\beta_1, \beta_1\gamma \rangle
\]
Then $A$ is a gentle algebra \cite{AS1987} and its marked surface $\bfS$ is shown in Figure \ref{fig:MS}.
Here, the marked surface is the geometric model which is used to describe the module category of a gentle algebra,
see for example \cite{BCS2021}.
\begin{figure}[htbp]
  \centering
\begin{tikzpicture}[scale=1.2]
\draw[line width = 1pt] (2,0) arc(0:360:2);
\draw[line width = 1pt] (0.5,0) arc(0:360:0.5);
\draw[blue, fill=blue] ( 0.00, 2.00) circle(0.7mm);
\draw[blue, fill=blue] ( 0.00,-0.50) circle(0.7mm);
\draw[blue, fill=blue] (-1.41,-1.41) circle(0.7mm);
\draw[blue, fill=blue] ( 1.41,-1.41) circle(0.7mm);
\draw[red, fill=white] (-1.41, 1.41) circle(0.7mm) [line width = 1pt];
\draw[red, fill=white] ( 1.41, 1.41) circle(0.7mm) [line width = 1pt];
\draw[red, fill=white] ( 0.00,-2.00) circle(0.7mm) [line width = 1pt];
\draw[red, fill=white] ( 0.00, 0.50) circle(0.7mm) [line width = 1pt];
\draw[blue] (0,2) to[out=-135,in=90] (-1,0.2) to[out=-90,in=-135] (0,-0.5);
\draw[blue] (0,2) to[out=-45, in=90] ( 1,0.2) to[out=-90, in=-45] (0,-0.5);
\draw[blue] (0,2) to[out=-160,in=90] (-1.4,0) -- (-1.4,-1.4);
\draw[blue] (0,2) to[out=-20, in=90] ( 1.4,0) -- ( 1.4,-1.4);
\draw[orange][line width = 1pt] (-1.4,1.4) to[out=20,in=160] ( 1.4,1.4);
\draw[orange][line width = 1pt] (0,0) circle(1cm);
\draw[orange][->] (0,1.68) -- (0.4,1.28) -- (2,1.28)
  node[right]{\tiny permissible curve with endpoint ``$\rmark$''};
\draw[orange][->] (0,-1.0) -- (-0.4,-1.2) -- (-2,-1.2)
  node[left]{\tiny permissible curve without endpoint};
\draw[orange]     (-2,-1.5)
  node[left]{\tiny (permitted closed curve)};
\end{tikzpicture}
  \caption{The marked surface $\bfS$ of $A$}
  \label{fig:MS}
\end{figure}
In this figure, ``$\rmark$'' and ``$\bmark$'' are called {\it $\rmark$-marked points} and {\it $\bmark$-marked points}, respectively, and the curves with endpoint ``$\bmark$'' are called {\it $\bmark$-arcs}.
The set of all $\bmark$-arcs is written as $\Delta$ which corresponds one-to-one with $\Q_0$.
Recall that a {\it permissible curve} on $\bfS$ is one of the following curve.
\begin{itemize}
  \item[(I)] A curve $c$ whose endpoints are $\rmark$-marked points such that
    each segment obtained by $\Delta$ cutting $c$ is an {\it arc segment},
    i.e., it is one of the forms (A) and (B) shown in Figure \ref{fig:arc segment}.

\begin{figure}[htbp]
  \centering
\begin{tikzpicture}
\draw[blue, fill=blue] (0,0) circle(0.7mm);
\draw[blue] ( 1.5,-2) -- (0,0) -- (-1.5,-2);
\draw[orange][line width=1pt] (-0.75,-1) -- ( 0.75,-1);
\draw[orange] (0,-1) node[below]{\tiny arc segment};
\draw (0,-2.5) node{(A)};
\end{tikzpicture}
\ \
\begin{tikzpicture}
\draw[blue, fill=blue] (0,0) circle(0.7mm) (0,-2) circle(0.7mm);
\draw[blue] (0,0) -- (0,-2);
\draw[black][line width=1pt] (0,0) to[out=180,in=90] (-1,-1) to[out=-90,in=180] (0,-2);
\draw[orange][line width=1pt] (-1,-1) -- (0,-1);
\draw[orange] (0,-1) node[below]{\tiny arc segment};
\draw[red, fill=white][line width=1pt] (-1,-1) circle(0.7mm);
\draw (0,-2.5) node{(B)};
\end{tikzpicture}
  \caption{Arc segments}
  \label{fig:arc segment}
\end{figure}

  \item[(II)] A curve $c$ without endpoints such that each segment of it is an arc segment of type (B) shown in Figure \ref{fig:arc segment}.
      In this case, $c$ is a closed curve.
\end{itemize}
See \cite[Definitions 3.1 and 3.5]{BCS2021}.
All indecomposable right $A$-modules can be described by using permissible curves.
To be more precise, we have the following bijection
\[ \X: \PC(\bfS) \cup (\CC(\bfS)\times\mathscr{J}) \to \ind(\modcat A), \]
where $\PC(\bfS)$ is the set of all permissible curve of type (A),
$\CC(\bfS)$ is the set of all closed curve of type (B),
and $\mathscr{J}$ is the set of all Jordan block with eigenvalue $\lambda\in\kk\backslash\{0\}$.
Each indecomposable modules lying in $\mathrm{Im}(\X|_{\PC(\bfS)})$ are called {\it string modules},
and each indecomposable modules lying in $\mathrm{Im}(\X|_{\CC(\bfS)\times\mathscr{J}})$ are called {\it band modules}.
Notice that the marked surface $\bfS$ can be used to describe the derived category of $A$,
see for example \cite[etc]{APS2023, HKK2017, LZ2021, OPS2018, QZZ2022}.

By the above correspondence $\X$, the curves $c_1$, $c_2$, $c_3$ and $c_4$ shown in Figure \ref{fig:stautilt} (1)
respectively describe the indecomposable right $A$-modules
$
\left(\begin{smallmatrix}
1 \\ 2\\ 3\\ 4
\end{smallmatrix}\right)_A, \
\left(\begin{smallmatrix}
1 \\ 2\\ 3
\end{smallmatrix}\right)_A, \
\left(\begin{smallmatrix}
1 \\ 2
\end{smallmatrix}\right)_A \ \text{and}\
(1)_A.
$
Then the curves $c_1$, $c_2$, $c_3$ and $c_4$ form a dissection $\Gamma$ of $\bfS$ which corresponding to the right $A$-module
\begin{align}\label{formula:stautilt}
  T = \left(\begin{smallmatrix}
1 \\ 2\\ 3\\ 4
\end{smallmatrix}\right)_A
\oplus \left(\begin{smallmatrix}
1 \\ 2\\ 3
\end{smallmatrix}\right)_A
\oplus \left(\begin{smallmatrix}
1 \\ 2
\end{smallmatrix}\right)_A
\oplus (1)_A.
\end{align}
By \cite[Theorem B]{HZZ2023}, we have that $T$ is a support $\tau$-tilting module
($\Gamma$ is said to be a generalized dissection in \cite{HZZ2023}).
Moreover, $A$ as a right $A$-module is also a support $\tau$-tilting module,
the dissection corresponding to it is shown in Figure \ref{fig:stautilt} (2).
\begin{figure}[htbp]
  \centering
\begin{tikzpicture}[scale=1.2]
\draw[line width = 1pt] (2,0) arc(0:360:2);
\draw[line width = 1pt] (0.5,0) arc(0:360:0.5);
\draw[orange][line width=1pt]
  (-1.4, 1.4) to[out=30,in=150] ( 1.4, 1.4);
\draw[orange][line width=1pt]
  (-1.4, 1.4) to[out=-10,in=170]
  ( 0.0, 1.1) to[out=-10,in=90]
  ( 1.0, 0.0) to[out=-90,in=90]
  ( 0.0,-2.0);
\draw[orange][line width=1pt]
  (-1.4, 1.4) -- ( 0.0, 0.5);
\draw[orange][line width=1pt]
  (-1.4, 1.4) to[out=-115,in=170] ( 0.0,-2.0);
\draw[orange]
  ( 0.0, 1.9) node[below]{$c_1$}
  ( 0.0, 1.0) node[above]{$c_2$}
  (-0.5, 0.9) node[below]{$c_3$}
  (-1.3,-1.1) node[right]{$c_4$}
;
\draw[blue, fill=blue] ( 0.00, 2.00) circle(0.7mm);
\draw[blue, fill=blue] ( 0.00,-0.50) circle(0.7mm);
\draw[blue, fill=blue] (-1.41,-1.41) circle(0.7mm);
\draw[blue, fill=blue] ( 1.41,-1.41) circle(0.7mm);
\draw[red, fill=white] (-1.41, 1.41) circle(0.7mm) [line width = 1pt];
\draw[red, fill=white] ( 1.41, 1.41) circle(0.7mm) [line width = 1pt];
\draw[red, fill=white] ( 0.00,-2.00) circle(0.7mm) [line width = 1pt];
\draw[red, fill=white] ( 0.00, 0.50) circle(0.7mm) [line width = 1pt];
\draw[blue] (0,2) to[out=-135,in=90] (-1,0.2) to[out=-90,in=-135] (0,-0.5);
\draw[blue] (0,2) to[out=-45, in=90] ( 1,0.2) to[out=-90, in=-45] (0,-0.5);
\draw[blue] (0,2) to[out=-160,in=90] (-1.4,0) -- (-1.4,-1.4);
\draw[blue] (0,2) to[out=-20, in=90] ( 1.4,0) -- ( 1.4,-1.4);
\draw (0,-2.5) node{(1)};
\end{tikzpicture}
\ \ \ \
\begin{tikzpicture}[scale=1.2]
\draw[line width = 1pt] (2,0) arc(0:360:2);
\draw[line width = 1pt] (0.5,0) arc(0:360:0.5);
\draw[orange][line width=1pt]
  (-1.4, 1.4) to[out=  30,in= 150] ( 1.4, 1.4);
\draw[orange][line width=1pt]
  ( 0.0, 0.5) to[out=  45,in=  90]
  ( 1.2, 0.0) to[out= -90,in=   0]
  ( 0.0,-1.2) to[out= 180,in= -90]
  (-1.2, 0.0) to[out=  90,in=-180]
  ( 0.0, 1.2) to[out=   0,in= -90]
  ( 1.4, 1.4);
\draw[orange][line width=1pt]
  ( 0.0, 0.5) to[out=  50,in= -65]
  ( 1.4, 1.4);
\draw[orange][line width=1pt]
  ( 0.0,-2.0) to[out=  10,in= -45]
  ( 1.4, 1.4);
\draw[blue, fill=blue] ( 0.00, 2.00) circle(0.7mm);
\draw[blue, fill=blue] ( 0.00,-0.50) circle(0.7mm);
\draw[blue, fill=blue] (-1.41,-1.41) circle(0.7mm);
\draw[blue, fill=blue] ( 1.41,-1.41) circle(0.7mm);
\draw[red, fill=white] (-1.41, 1.41) circle(0.7mm) [line width = 1pt];
\draw[red, fill=white] ( 1.41, 1.41) circle(0.7mm) [line width = 1pt];
\draw[red, fill=white] ( 0.00,-2.00) circle(0.7mm) [line width = 1pt];
\draw[red, fill=white] ( 0.00, 0.50) circle(0.7mm) [line width = 1pt];
\draw[blue] (0,2) to[out=-135,in=90] (-1,0.2) to[out=-90,in=-135] (0,-0.5);
\draw[blue] (0,2) to[out=-45, in=90] ( 1,0.2) to[out=-90, in=-45] (0,-0.5);
\draw[blue] (0,2) to[out=-160,in=90] (-1.4,0) -- (-1.4,-1.4);
\draw[blue] (0,2) to[out=-20, in=90] ( 1.4,0) -- ( 1.4,-1.4);
\draw (0,-2.5) node{(2)};
\end{tikzpicture}
  \caption{A support $\tau$-tilting module}
  \label{fig:stautilt}
\end{figure}

In \cite{WW1985}, Wald and Waschb\"{u}sch introduced V-sequences and primitive V-sequences to describe the finitely generated module category of a biserial algebra,
and provided a classification of all indecomposable modules.
Furthermore, if a biserial algebra is a string algebra, then Butler and Ringel showed that the above description is a bijection \cite{BR1987}.
In this case, V-sequences and primitive V-sequences are called strings and bands.
In derived category, the concepts similar to strings and bands respectively are homotopy strings and homotopy bands,
which are used to study derived representation types of gentle algebras, see for example \cite[etc]{CZ2019, Z2019, ZH2016}.
Now, we recall strings and bands on a bound quiver $(\Q,\I)$.

For any arrow $\alpha$, its {\it formal inverse} is written as $\alpha^{-1}$,
and, naturally, we can define $\target(\alpha^{-1}) := \source(\alpha)$, $\target(\alpha^{-1})=\source(\alpha)$,
$\Q_1^{-1} := \{\alpha^{-1} \mid \alpha \in \Q_1\}$, and $\Q_1^{\pm 1} := \Q_1\cup\Q_1^{-1}$.

A {\it string} on the bound quiver of a gentle algebra is a sequence $s=a_1a_2\cdots a_n$ ($n\in\mathbb{N}$) of elements lying in $\Q_1^{\pm 1}$ such that:
\begin{itemize}
  \item $\target(a_i)=\source(a_{i+1})$ holds for all $1\leqslant i\leqslant n-1$;
  \item if $a_i, a_{i+1} \in \Q_1$, then $a_ia_{i+1}\notin\I$, and if $a_i, a_{i+1} \in \Q_1^{-1}$, then $a_{i+1}^{-1}a_{i}^{-1}\notin\I$;
  \item if $a_i\in\Q_1$ and $a_{i+1}\in\Q_1^{-1}$, then $a_i\ne a_{i+1}^{-1}$,
    and if $a_i\in\Q_1^{-1}$ and $a_{i+1}\in\Q_1$, then $a_i^{-1}\ne a_{i+1}$.
\end{itemize}
In particular, in the case for $n=0$, $s$ is a path of length zero which is seen as a string of length zero.
Moreover, we define the empty $\varnothing$ is a {\it trivial string} which is used to correspond the zero module.
Obviously, the formal inverse $s^{-1}=a_n^{-1}\cdots a_2^{-1}a_1{-1}$ is also a string.
We call that two string $s$ and $s'$ are {\it equivalent} if $s=s'$ or $s^{-1}=s'$.

A {\it band} $b=b_0b_1\cdots b_{n-1}$ is a string such that:
\begin{itemize}
  \item $\target(b_{n-1}) = \source(b_0)$;
  \item $b^2$ is a string, that is, $b_{n-1}b_1\notin\I$;
  \item $b$ is not a non-power of any string, that is, $b\ne s^n$ for any string $s$ and any $n\geqslant 2$.
\end{itemize}
We call that tow bands $b$ and $b'$ are {\it equivalent} if there exists $0\leqslant t\leqslant n-1$ such that one of  $b[t]=b'$ and $(b[t])^{-1}=b'$ holds.
Here, $b[t]$ is the band $b_{t}b_{t+1}\cdots b_{n-1}b_0b_1\cdots b_{t-1}$.

Let $\Str(A)$ be the set of all equivalent classes of strings on $(\Q_A,\I_A)$
and $\Band(A)$ be the set of all equivalent classes of bands on $(\Q_A,\I_A)$.
We have the following bijection:
\[ \M: \Str(A) \cup (\Band(A)\times\mathscr{J}) \to \ind(\modcat A). \]
Furthermore, we have $\Im(\X|_{\PC(\bfS)}) = \Im(\M|_{\Str(A)})$ and
$\Im(\X|_{\CC(\bfS)\times\mathscr{J}}) = \Im(\M|_{\Band(A)\times\mathscr{J}})$,
c.f. \cite[Section 3]{BR1987} and \cite[Theorems 3.8 and 3.9]{BCS2021}.

Next, we provide an instance for Theorem \ref{main1}. In our examples, all modules are basic,
that is, for any two indecomposable direct summand $N_1$ and $N_2$ of any module,
we have $N_1\not\cong N_2$.

\begin{example} \label{examp:for mian1}
Take $e=e_1+e_2$, where for any $i\in (\Q_A)_0$, $e_i$ is the primitive idempotent corresponding to the vertex $i$.
Then we a recollement by \cite[Example 2.7]{P} as follows.
$$
\xymatrix{
    \modcat B \ar[rr]^{i_* = \mathrm{inc}}
& & \modcat A
    \ar@/_1.5pc/[ll]_{i^*=-\otimes_{A}B}
    \ar@/^1.5pc/[ll]^{i^!=\mathrm{Hom}_{A}(B,-)}
    \ar[rr]^{j^*=(-)e}
& & \modcat C
    \ar@/_1.5pc/[ll]_{j_!=-\otimes_{C}eA}
    \ar@/^1.5pc/[ll]^{j_*=\mathrm{Hom}_{C}(Ae,-)}
    }
$$
Here, $B = A/AeA$ is isomorphic to the bound quiver algebra $\kk\Q_B/\I_B$ whose bound quiver is
\[ \Q_B = \ \xymatrix{3 \ar[r]^{\gamma} & 4 }, \ \ \I_B = 0, \]
and $C = eAe$ is isomorphic to the bound quiver algebra $\kk\Q_C/\I_C$ whose bound quiver is
\[ \Q_C = \ \xymatrix{1 \ar[r]^{\alpha} & 2 }, \ \ \I_C = 0. \]

We have the isomorphism of $\kk$-linear spaces $B = A/AeA \cong \kk e_3 + \kk e_4 + \kk \gamma$.
Thus, the quiver representation of $B$, as a right $A$-module, is
\[
\bigg(\xymatrix{
  Be_1 \ar[r]^{\alpha}
& Be_2 \ar@/^0.5pc/[r]^{\beta_1}
    \ar@/_0.5pc/[r]_{\beta_2}
& Be_3 \ar[r]^{\gamma}
& Be_4
}\bigg)
\cong
\bigg(\xymatrix{
  0 \ar[r]
& 0 \ar@/^0.5pc/[r]
    \ar@/_0.5pc/[r]
& \kk \ar[r]^{[{^1_0}]}
& \kk^{\oplus 2}
}\bigg),
\]
and so, $B$ is isomorphic to the indecomposable projective right $A$-module $P(3)_A\oplus P(4)_A$.
Then $i^* = -\otimes_A B$ and $i^!=\Hom_A(B,-)$ are exact.

For the algebra $B$, we have $\add(({^3_4})_B \oplus (4)_B)$ is Wakamatsu tilting,
and for the algebra $C$, we have $\add (({^1_2})_C \oplus (2)_C)$ is Wakamatsu tilting.
Let $X \in \modcat A$ such that
\begin{align}\label{formula:examp 1}
  i^*(X) \in \add(({^3_4})_B \oplus (4)_B)
\end{align}
and
\begin{align}\label{formula:examp 2}
  j^*(X) \in \add (({^1_2})_C \oplus (2)_C).
\end{align}
Then we obtain that any indecomposable direct summand of $X$ is isomorphic to one of the following forms:
\begin{align}\label{formula:examp 3}
\M(e_4),\ 
\M(\gamma),\ 
\M((\beta_2\beta_1^{-1})^{t_1}\beta_2\gamma),\ 
\M(\beta_1^{-1}(\beta_2\beta_1^{-1})^{t_2}\beta_2\gamma),\ 
\text{and}\
\M(\alpha(\beta_2\beta_1^{-1})^{t_3}\beta_2\gamma),\ 
\end{align}

Let $e_{B,3} = e_3+AeA$ and $e_{B,4} = e_4+AeA$ be two idempotent corresponding to the vertices $3,4\in (\Q_B)_0$. Then we have:
\begin{align*}
    i^*(\M((\beta_2\beta_1^{-1})^{t_1}\beta_2\gamma)) e_{B,3}
& = \M((\beta_2\beta_1^{-1})^{t_1}\beta_2\gamma) \otimes_A Be_{B,3} \\
& = \M((\beta_2\beta_1^{-1})^{t_1}\beta_2\gamma) \otimes_A (\kk e_{3}+AeA) \\
& \cong \M((\beta_2\beta_1^{-1})^{t_1}\beta_2\gamma) \otimes_{\kk} \kk e_{3} + AeA \\
& \cong \kk^{\oplus (t_1+1)}.
\end{align*}
Since $X$ is basic, we obtain $t_1=0$, and so $\M((\beta_2\beta_1^{-1})^{t_1}\beta_2\gamma)
= \M(\beta_2\gamma) \cong P(3)_A$.
One can check that $t_2=0$ and $t_3=0$ by the similar way.
Thus, each indecomposable direct summand of $X$ is isomorphic to one of
\begin{align}\label{formula:examp 4}
\M(e_4) \cong (4)_A,\ 
\M(\gamma) \cong ({^3_4})_A,\ 
\M(\beta_1^{-1}\beta_2\gamma)
  \cong
  \left(
    \begin{smallmatrix}
     & 2 & \\ 3 & & 3 \\  & & 4
    \end{smallmatrix}
  \right)_A, \ 
\text{and}\
\M(\alpha\beta_2\gamma)
  \cong
  \left(
    \begin{smallmatrix}
    1 \\ 2\\ 3\\ 4
    \end{smallmatrix}
  \right)_A.\ 
\end{align}
Since $j^*(\M(e_4)) = 0$, $j^*(\M(\gamma)) = 0$,
$j^*(\M(\beta_1^{-1}\beta_2\gamma) ) \cong S(2)_C$,
$j^*(\M(\alpha\beta_2\gamma) ) \cong ({^1_2})_C$
are objects in $\add (({^1_2})_C \oplus (2)_C)$,
we obtain that
\[ X = \bigoplus_{M\in(\ref{formula:examp 4})} M \]
and $\add X$ is a Wakamatsu tilting subcategory of $\modcat A$ by Theorem \ref{main1}. \
Indeed, $X$ is also a support $\tau$-tilting module (see Figure \ref{fig:stautilt} (2)) in this example, and it is isomorphic to $A_A$.
\end{example}

The recollement given in Example \ref{examp:for mian1} is also an instance for Theorem \ref{main2},
see the following example.

\begin{example} \label{examp:for mian2}
All indecomposable $A$-modules can be divided to five classes as follows by the correspondence $\M$:
\begin{itemize}
  \item[(1)] simple modules: $S(1)_A$, $S(2)_A$, $S(3)_A$ and $S(4)_A$;
  \item[(2)] $\M(\alpha(\beta_2\beta_1^{-1})^t)$, $\M(\alpha(\beta_2\beta_1^{-1})^t\beta_2)$
    and $\M(\alpha(\beta_2\beta_1^{-1})^t\beta_2\gamma)$;
  \item[(3)] $\M((\beta_2\beta_1^{-1})^t)$, $\M((\beta_2\beta_1^{-1})^t\beta_2)$
    and $\M((\beta_2\beta_1^{-1})^t\beta_2\gamma)$;
  \item[(4)] $\M(\gamma)$;
  \item[(5)] band modules: $\M(\beta_2\beta^{-1}, \pmb{J}_n(\lambda))$,
    where $\pmb{J}_n(\lambda)) =
    \left(
      \begin{smallmatrix}
       \lambda & 1 & & \\
       & \lambda & \cdots & \\
       & & \cdots & 1 \\
       & & &  \lambda
      \end{smallmatrix}
    \right)$, $\lambda\ne 0$.
\end{itemize}
The module $X = \bigoplus\limits_{M\in(\ref{formula:examp 4})} M $ given in Example \ref{examp:for mian1}
decides a subcategory $\add X$ of $\modcat A$ which is Wakamatsu tilting.

Next, for any indecomposable module $N$ in ${^{\bot}(\add X)}$,
we show $i_*i^!(N) \in \add X$.
We will consider all $i_*i^!(N)$ for all indecomposable right $A$-modules $N$ in this proof.
Then, obviously, if $i_*i^!(N)=0$ or $i_*i^!(N) \in \add X$,
then we don't need to worry about whether $N \in {^{\bot}(\add X)}$ holds true or not.

\begin{itemize}
\item[(1)]
For the indecomposable module in the case (1):
\begin{itemize}
  \item[(1.1)] $i_*i^! (S(1)_A) = i_*\Hom_A(A/AeA, S(1)_A)=0$.
  \item[(1.2)] $i_*i^! (S(2)_A) = i_*\Hom_A(A/AeA, S(2)_A)=0$.
  \item[(1.3)] $i_*i^! (S(3)_A) = i_*\Hom_A(A/AeA, S(3)_A)\ne 0$.
    However, we have $S(3)_A \notin {^{\bot}(\add X)}$
    since there is a non-split short exact sequence $0 \to (4)_A \to ({^3_4})_A \to (3)_A \to 0$
    with $(4)_A \in (\ref{formula:examp 4})$.
    Thus, this case is not something we should consider.

  \item[(1.4)] $i_*i^! (S(4)_A) = i_*\Hom_A(A/AeA, S(4))\ne 0$.
    One can check that $i_*\Hom_A(A/AeA, S(4)) \cong S(4)_A$ by the following way:
    \begin{align*}
          & \Hom_A(A/AeA, S(4)_A) e_4 \\
       =\ & \Hom_A(\kk e_3 + \kk e_4 + \kk \gamma + AeA, \kk e_4 ) e_4 \\
   \cong\ & \Big((k_{e_4} e_4 + AeA) \mapsto (k_{e_4} e_4 ) \Big)e_4 ~~ (k_{e_4}\in \kk) \\
   \cong\ & S(4)_\kk.
    \end{align*}
    By $S(4)_A \in (\ref{formula:examp 4})$, we have $i_*i^! (S(4)_A)  \cong S(4)_A \in \add X$ as required.
\end{itemize}

\item[(2)]
For the indecomposable module in the case (2):
\begin{itemize}
  \item[(2.1)] We have a non-split short exact sequence
    \[0 \to ({^3_4})_A \to \M(\alpha(\beta_2\beta_1^{-1})^t\beta_2\gamma)
        \to \M(\alpha(\beta_2\beta_1^{-1})^t) \to 0\]
    for each $t\geqslant 0$. Here, we have $ ({^3_4})_A \cong \M(\gamma) \in (\ref{formula:examp 4})$,
    and obtain $\M(\alpha(\beta_2\beta_1^{-1})^t) \notin {^{\bot}(\add X )}$ as required.

  \item[(2.2)] We have a non-split short exact sequence
    \[0 \to (4)_A \to \M(\alpha(\beta_2\beta_1^{-1})^t\beta_2\gamma)
        \to \M(\alpha(\beta_2\beta_1^{-1})^t\beta_2) \to 0\]
    for each $t\geqslant 0$. Here, we have $(4)_A \cong \M(e_4) \in (\ref{formula:examp 4})$,
    and obtain $\M(\alpha(\beta_2\beta_1^{-1})^t\beta_2) \notin {^{\bot}(\add X )}$ as required.

  \item[(2.3)] For the string module $\M(\alpha(\beta_2\beta_1^{-1})^t\beta_2\gamma)
   = \left( \begin{smallmatrix}
    1 &   &   &   &   &   &   & \\
    2 &   & 2 &...&...& 2 &   & \\
      & 3 &   & 3 &...&...& 3 & \\
      &   &   &   &   &   & 4 &
   \end{smallmatrix} \right)_A$,
    we consider its projective resolution is
      \[ 0
      \to P(3)^{\oplus (t-1)}
      \to P(1) \oplus P(2)^{\oplus t}
      \to \M(\alpha(\beta_2\beta_1^{-1})^t\beta_2\gamma) \to 0. \]
    Let $\pmb{P}(\M(\alpha(\beta_2\beta_1^{-1})^t\beta_2\gamma))$ be the delete complex corresponding to the above resolution, then
    \[\Hom_{D^b(A)}
      (
        \pmb{P}(\M(\alpha(\beta_2\beta_1^{-1})^t\beta_2\gamma)),
        P(1)[1]
      ) \ne 0,\]
    where $D^b(A)$ is the derived category of $A$.
    It follows that $\M(\alpha(\beta_2\beta_1^{-1})^t\beta_2\gamma) \notin {^{\bot}P(1)}$,
    then we have $\M(\alpha(\beta_2\beta_1^{-1})^t\beta_2\gamma) \notin {^{\bot}(\add X)}$.
\end{itemize}

\item[(3)]
We can show that all indecomposable modules in (3) do not lie in ${^{\bot}(\add X)}$ by the way similar to (2).

\item[(4)]
For the indecomposable module in the case (4), it is isomorphic to $(^3_4)_A \cong P(3)_A = e_3A$.
Then, by $B\cong P(3)_A\oplus P(4)_A$, we have
\begin{align*}
  i_*i^!((^3_4)_A) \cong \Hom_A(e_3A \oplus e_4A , e_3A) \cong e_3Ae_3 \oplus_{\kk} e_3Ae_4 \cong P(3)_A \in \add X
\end{align*}
as required (``$\oplus_{\kk}$'' is the direct sum of $\kk$-linear space).

\item[(5)]
For each band module $\M(\beta_2\beta^{-1}, \pmb{J}_n(\lambda)) =: B(n,\lambda)$ in the case (5),
we have the projective dimension $\mathrm{proj.dim} B(n,\lambda) = 1$
and its $1^{\mathrm{st}}$-syzygy $\Omega^1(B(n,\lambda)) \cong P(2)^{\oplus n}$ is projective.
It follows that $\Hom_{D^b(A)}(\pmb{P}(B(n,\lambda)), P(v)[1])\ne 0$ holds for all $v\in \{2,3,4\}$.
Then $B(n,\lambda) \notin {^{\bot}(\add X)}$ as required.
\end{itemize}

Therefore, we have $i_*i^!({^{\bot}(\add X)}) \subseteq \add X$ in this example.
One can check that $j_*j^*(\add X) \subseteq \add X$.
Thus, by Theorem \ref{main2}, we obtain that
$i^!(\add X) = \add (P(3)_B\oplus P(4)_B)$ and $j^*(\add X) = \add (P(1)_C\oplus P(2)_C)$ are Wakamatsu tilting.
\end{example}

The following example provide an instance for Theorem \ref{main3}.

\begin{example} \label{examp:for mian3}
Consider the recollement given in Example \ref{examp:for mian1}.
Take the support $\tau$-tilting module $(3)_B \oplus ({^3_4})_B$ in $\modcat B$
and the support $\tau$-tilting module $(2)_C \oplus ({^1_2})_C$ in $\modcat C$.
Then $\mathcal{Z}' := \add ((3)_B \oplus ({^3_4})_B)$ and
$\mathcal{Z}'' := \add ((1)_C \oplus ({^1_2})_C)$ are
weak support $\tau$-tilting subcategories of $\modcat B$ and $\modcat C$, respectively.
By using the method similar to Example \ref{examp:for mian1}, we have
\begin{align*}
    \mathcal{Z}
& = \bigg\{
      Z \in \modcat B
        \mid Z\otimes_A B \in \mathcal{Z}' \text{~and~}
             Ze \in \mathcal{Z}''
    \bigg\} \\
& = \add\bigg(
    \left(\begin{smallmatrix}
     1 \\ 2 \\ 3 \\ 4
    \end{smallmatrix}\right)_A
    \oplus
    \left(\begin{smallmatrix}
     1 \\ 2 \\ 3
    \end{smallmatrix}\right)_A
    \oplus
    \left(\begin{smallmatrix}
     1 \\ 2
    \end{smallmatrix}\right)_A
    \oplus
    (1)_A
    \bigg) \\
& \mathop{=}\limits^{(\ref{formula:stautilt})} \add T
\end{align*}
Thus, $\mathcal{Z}$ is a weak support $\tau$-tilting subcategory of $\modcat A$ by Theorem \ref{main3}.
\end{example}

Now we provide the last example as follows which is an instance for Theorem \ref{main4}.

\begin{example} \label{examp:for mian4}
Consider the recollement given in Example \ref{examp:for mian1}.
Take the subcategory $\add T$ with $T = (\ref{formula:stautilt})$.
Then it is a weak support $\tau$-tilting subcategory of $\modcat A$.
Moreover, by the method similar to Example \ref{examp:for mian2},
one can check that $i_*i^*(\add T)\subseteq \add T$
and $j_*j^*(\add T)\subseteq \add T$.
Then
\[ i^*(\add T) = \add (i^*(T)) = \add((1)_B \oplus ({^1_2})_B) \]
and
\[ j^*(\add T) = \add (j^*(T)) = \add((3)_C \oplus ({^3_4})_C) \]
are weak support $\tau$-tilting subcategories of $\modcat B$ and $\modcat C$, respectively.
\end{example}

\vspace{0.5cm}

\noindent
\textbf{Funding}
\hspace{2mm}
Yongduo Wang is supported by National Natural Science Foundation of China (Grant No. 11861043).
Yu-Zhe Liu was supported by the National Natural Science Foundation of China (Grant Nos. 12401042 and 12171207),
Guizhou Provincial Basic Research Program (Natural Science) (Grant Nos. ZK[2024]YiBan066 and ZD[2025]085),
and Scientic Research Foundation of Guizhou University (Grant Nos. [2022]53, [2022]65, [2023]16).
Jian He is supported by the National Natural Science Foundation of China (Grant No. 12501048) and Youth
Science and Technology Foundation of Gansu Provincial (Grant No. 23JRRA825).
Dejun Wu is supported by National Natural Science Foundation of China (Grant No. 12261056).
\vspace{3mm}

\noindent
\textbf{Data Availability}
\hspace{2mm}
Data sharing not applicable to this article as no datasets were generated or analysed during the current study.
\vspace{3mm}

\noindent
\textbf{Conflict of Interests}
\hspace{2mm}
The authors declare that they have no conflicts of interest to this work.
\vspace{3mm}

\noindent
\textbf{Acknowledgements}
The authors would like to express their gratitude to the referee for providing valuable comments and suggestions that have helped shape the article into its current form.

\hspace{2mm}

\textbf{Yongduo Wang}\\
1. School of Mathematics and Statistics, Gansu Normal University for Nationalities, 747000 Hezuo, Gansu, P. R. China\\
2. Department of Applied Mathematics, Lanzhou University of Technology, 730050 Lanzhou, Gansu, P. R. China\\
E-mail: \textsf{ydwang@lut.edu.cn}\\[0.3cm]
\textbf{Hongyang Luo}\\
Department of Applied Mathematics, Lanzhou University of Technology, 730050 Lanzhou, Gansu, P. R. China\\
E-mail: \textsf{ 232070101001@lut.edu.cn}\\[0.3cm]
\textbf{Yu-Zhe Liu}\\
School of Mathematics and statistics, Guizhou University, 550025 Guiyang, Guizhou, P. R. China\\
E-mail: \textsf{liuyz@gzu.edu.cn / yzliu3@163.com}\\[0.3cm]
\textbf{Jian He}\\
Department of Applied Mathematics, Lanzhou University of Technology, 730050 Lanzhou, Gansu, P. R. China\\
E-mail: \textsf{jianhe30@163.com}\\[0.3cm]
\textbf{Dejun Wu}\\
Department of Applied Mathematics, Lanzhou University of Technology, 730050 Lanzhou, Gansu, P. R. China\\
E-mail: \textsf{ wudj@lut.edu.cn}

\end{document}